\documentclass[a4paper,11pt,makeidx]{amsart}
\usepackage{amscd}
\usepackage{xypic}  %commu 
\usepackage{amssymb}
\usepackage{amsthm}
\usepackage{epsf}
\makeindex
\newtheoremstyle{break}% name
  {9pt}%      Space above, empty = `usual value'
  {9pt}%      Space below
  {\itshape}% Body font
  {}%         Indent amount (empty = no indent, \parindent = para indent)
  {\bfseries}% Thm head font
  {.}%        Punctuation after thm head
  {\newline}% Space after thm head: \newline = linebreak
  {}%         Thm head spec
\theoremstyle{break}
\newtheorem{bthm}{Theorem}
\newtheorem{bprop}{Proposition}
\newtheorem{bcor}{Corollary}
\theoremstyle{plain}
\newtheorem{thm}{Theorem}[section]

\newtheorem{lemma}[thm]{Lemma}
\newtheorem{sublemma}[thm]{Sublemma}
\newtheorem{prop}[thm]{Proposition}

\newtheorem{defn}[thm]{Definition}
\newtheorem{claim}[thm]{Claim}
\newtheorem{rem}[thm]{Remark}

\renewcommand{\proofname}{Proof}

\def\Hilb{\operatorname{Hilb}}

\def\min{\operatorname{min}}
\def\Im{\operatorname{Im}}

\def\cork{\operatorname{cork}}
\def\max{\operatorname{max}}
\def\length{\operatorname{length}}
\def\c1{\operatorname{c_1}}
\def\c2{\operatorname{c_2}}
\def\Cliff{\operatorname{Cliff}}

\def\gon{\operatorname{gon}}
\def\gengon{\operatorname{gengon}}
\def\mingon{\operatorname{mingon}}

\def\PP{{\mathbb P}}

\def\N{{\mathcal N}}
\def\O{{\mathcal O}}
\def\I{{\mathcal J}}

\def\E{{\mathcal E}}

\def\U{{\mathcal U}}
\def\V{{\mathcal V}}
\def\J{{\mathfrak J}}
\def\*{\otimes}
\def\eqv{\equiv}
\def\sub{\subseteq}

\def\+{\oplus}                   % direct sum
\def\*{\otimes}                  % tensor product
\def\hpil{\longrightarrow}       % ----->

\def\Pic{\operatorname{Pic}}

\def\det{\operatorname{det}}
\def\Bs{\operatorname{Bs}}

\begin{document}

\title[Brill-Noether theory of curves on Enriques surfaces I]{Brill-Noether theory of curves on
Enriques surfaces I: \\ the positive cone and gonality} 

\dedicatory{\normalsize \dag \ Dedicated to the memory of Silvano Bispuri} 

\author[A.L. Knutsen and A.F. Lopez]{Andreas Leopold Knutsen* and Angelo Felice Lopez**}

\address{\hskip -.43cm Andreas Leopold Knutsen, Department of Mathematics, University of Bergen,
Johannes Brunsgate 12, 5008 Bergen, Norway. e-mail {\tt andreas.knutsen@math.uib.no}}

\thanks{* Research partially supported by a Marie Curie Intra-European Fellowship within the 6th 
European Community Framework Programme}

\address{\hskip -.43cm Angelo Felice Lopez, Dipartimento di Matematica, Universit\`a di Roma 
Tre, Largo San Leonardo Murialdo 1, 00146, Roma, Italy. e-mail {\tt lopez@mat.uniroma3.it}}

\thanks{** Research partially supported by the MIUR national project ``Geometria delle variet\`a 
algebriche" COFIN 2002-2004.}

\thanks{{\it 2000 Mathematics Subject Classification} : Primary 14H51, 14C20, 14J28. Secondary 14J05,
14F17}

\begin{abstract} 

We study the existence of linear series on curves lying on an Enriques surface and general in
their complete linear system. Using a method that works also below the Bogomolov-Reider
range, we compute, in all cases, the gonality of such curves. We also give a new result 
about the positive cone of line bundles on an Enriques surface and we show how this relates to
the gonality.   
\end{abstract}

\maketitle

\section{Introduction}
\label{intro}

Let $S$ be a smooth surface and let $L$ be a line bundle on $S$. There is a natural 
interaction between the geometry of $S$ and that of the curves $C \in |L|$. On the one hand, strong geometric
properties of curves in $|L|$, do lead, in many cases, to analogous properties of $S$ itself. This is the case
for example if all smooth curves in $|L|$ are hyperelliptic (\cite{ca, en, sv}), or trigonal
(\cite{se, re, pao, fa}), etc.. On the other hand one can choose an
interesting type of surface and try to derive as much information as possible on the
curves in $|L|$. Celebrated examples of this investigation range from well-known classical ones
to very recent ones. Perhaps this line of thought was revived by the Green-Harris-Mumford
conjecture, namely that all smooth curves in a given linear series on a K3 surface have the same
Clifford index.

To study curves on K3 surfaces new interesting vector bundle methods were introduced by
Lazarsfeld, Tyurin, Reider, Donagi and Morrison \cite{dm}, culminating on one side with the proof
of the Green-Harris-Mumford conjecture \cite{gl} and on the other with the fact that curves on a
K3 surface of rank one do behave, from the point of view of Brill-Noether theory, like general
ones \cite{la}.

The study of Brill-Noether theory of curves in a given complete linear system $|L|$ on a
surface $S$ has, besides its own beauty, lots of interesting applications. We mention here the one that was
our main motivation. 

Suppose that $L$ is very ample, giving an embedding $S \subset \PP^r = \PP H^0(L)$. In
the study of threefolds it is interesting to know whether there exists a threefold
$X \subset \PP^{r + 1}$ different for a cone over $S$ and such that $S = X \cap \PP^r$. If $r
\geq 4$ there is a well-known condition (Zak's theorem \cite{za}): If
$h^0(N_{S/ \PP^r}(-1)) \leq r + 1$, where $N_{S/\PP^r}$ is the normal bundle of $S$, then there is
no such $X$. Now the cohomology of the normal bundle of $S$ is often related to the one of a smooth
hyperplane section $Y = S \cap H$. But on a curve we have the formula (\cite{wa})  
$h^0(N_{Y/ \PP^{r-1}}(-1)) = r + \cork \Phi_{H_Y, \omega_Y}$, where $\Phi_{H_Y, \omega_Y}$ is the
Gaussian map associated to the canonical and hyperplane bundle $H_Y$ of $Y$. At last, the
surjectivity of Gaussian maps on a curve $Y$ is very much governed by its Brill-Noether theory
(\cite{wa, bel}). Moreover, as the results of \cite{klm} show, this knowledge will be needed
also when $Y$ is not necessarily the hyperplane section of $X$.

In the present article we investigate the Brill-Noether theory of curves on an Enriques surface, emphasizing the
calculation of the gonality (see Thm. \ref{main}), prove a new result on the positive cone on an Enriques
surface (see Prop. \ref{poscone}) and show how these two results are related. In a subsequent paper \cite{kl2},
we will study Clifford index and exceptional curves. These results will play a crucial role in \cite{klm},
where, among other things, we prove the genus bound $g \leq 17$ for threefolds having an Enriques surface as
hyperplane section.

\vskip .2cm

Let now $S$ be an Enriques surface and let $|L|$ be a base-component free complete linear
system on $S$. Unlike (most cases of) $K3$ and Del Pezzo surfaces, it is not the case that the Clifford index or
the gonality of smooth curves in $|L|$ are constant, as simple examples show. Therefore perhaps the best first
question one can ask is about the linear series on {\it general} curves in $|L|$. Now recall from \cite{cd}:

\begin{defn}
Let $L$ be a line bundle on an Enriques surface $S$ such that $L^2 > 0$. Set
\[ \phi(L) := \inf \{|F.L| \; : \; F \in \Pic S, F^2 = 0, F \not\eqv 0\}. \] 
\end{defn} 

An important property of this function, which will be used throughout the article, is  
that $\phi(L)^2 \leq L^2$ \cite[Cor.2.7.1]{cd}. Hence, for $L^2>>0$, we have, for any smooth 
$C \in |L|$, that $\gon C \leq 2 \phi(L) << \lfloor \frac{L^2}{4} \rfloor + 2=
\lfloor \frac{g(C)+3}{2} \rfloor$, so that the curves are far from being Brill-Noether general, again unlike 
on (general) $K3$ surfaces. One might expect that either the Brill-Noether theory of general curves or
elliptic pencils on the surface are enough to calculate the gonality of general curves in $|L|$, but 
this turns out to be false: Let $|2E_1|, |2E_2|$ be two genus one pencils such that $E_1.E_2 = 2$ (they exist
on a general Enriques surface by \cite[Thm.2.7.2]{co1}) and let $L = n (E_1 + E_2)$. We have $L^2 = 4n^2$ and
$\phi(L) = 2n$. Let $C \in |L|$ be a general curve and set $B = E_1 + E_2$. For $n \geq 2$ we find that
$|B_{|C}|$ is a $g^2_{4n}$ that cannot be very ample, else $4n^2 = 4n(4n -3)$, hence $n = 1$. Therefore 
$\gon(C) \leq 4n  - 2 < 2\phi(L)$. 

Motivated by the above examples we give the ensuing

\begin{defn}
Let $L$ be a line bundle on an Enriques surface $S$ such that $L^2 > 0$. We define
\[ \mu(L) = \min\{B.L - 2 \ : \ B \in \Pic(S) \ \hbox{with B effective,} \ B^2 = 4, \phi(B) = 2, B
\not\eqv  L \}. \]
\end{defn} 

The first of the two main results of this paper shows that in fact the gonality of a
general smooth curve in a given complete linear system on an Enriques surface is governed by elliptic pencils
and divisors of self-intersection $4$:

\vskip .3cm

\begin{bthm}
\label{main}
Let $|L|$ be a base-component free complete linear system on an Enriques surface $S$ such that $L^2
> 0$. Then, for a general $C \in |L|$, we have
\[\gon(C) = \min\{2 \phi(L), \mu(L), \lfloor \frac{L^2}{4} \rfloor + 2 \}. \]
\end{bthm}

The question now arises of how to compute the ``new" function $\mu(L)$. In section 2 we
classify line bundles for which $\mu(L) < 2\phi(L)$ (Proposition \ref{prop:mu}). It turns out that such line
bundles are the ``extremal ones'' in the positive cone in a way we now explain. 

As mentioned above, $L^2 \geq \phi(L)^2$ for any line bundle $L$ on an Enriques surface with $L^2 > 0$. We prove
that there are no line bundles with $\phi(L)^2 < L^2 < \phi(L)^2 + \phi(L) - 2$ and we classify the borderline
cases. (The proposition is stated for simplicity for $L$ effective, otherwise it will hold for $K_S - L$ by
Riemann-Roch)

\vskip .2cm

\begin{bprop}
\label{poscone}
Let $L$ be an effective line bundle on an Enriques surface with $L^2 > 0$. If $L^2 \leq
\phi(L)^2 + \phi(L) - 2$, then there exist primitive effective divisors $E_i$ with $E_i^2 =
0$, for $i = 1, 2, 3$, $E_1.E_2 = E_1.E_3 = 2$, $E_2.E_3 = 1$ and an integer $h \geq 1$ so that one of
the two following occurs:
\begin{itemize}
\item[(i)] $L^2 = \phi(L)^2$. In this case $L \eqv h(E_1 + E_2)$.
\item[(ii)] $L^2 = \phi(L)^2 + \phi(L) - 2$. In this case either 
\begin{itemize}
\item[(ii-a)] $L \sim h(E_1 + E_2) + E_3$; or 
\item[(ii-b)] $L \sim (h + 1)E_1 + hE_2 + E_3$; or
\item[(ii-c)] $L \eqv 2(E_1 + E_2 + E_3)$ (whence $L^2 = 40$ and $\phi(L) = 6$).
\end{itemize}
\end{itemize}
\end{bprop}

The relation to the function $\mu(L)$ is that the line bundles appearing in (i), (ii-a) and
(ii-b) are precisely the ones for which $\mu(L) < 2\phi(L)$ (Proposition \ref{prop:mu}). In other words, linear
systems whose general members have ``nongeneral'' (in the sense of Brill-Noether theory) gonalities not
computed by elliptic pencils are precisely the ``extremal'' cases  (i), (ii-a) and (ii-b) in the positive
cone. 

In light of this, Theorem \ref{main} can be stated only in terms of $L^2$ and $\phi(L)$:

\vskip .2cm

\begin{bcor}
\label{main2}
Let $|L|$ be a base-component free complete linear system on an Enriques surface such that $L^2 > 0$
and  let $C \in |L|$ be a general curve. Then 
\[ \gon(C) = 2\phi(L) \]
unless $L$ is of one of the following types:
\begin{itemize}
\item[(a)] $L^2 = \phi(L)^2$ with $\phi(L) \geq 2$ and even. In these cases $\gon(C) = 2\phi(L) - 2$.
\item[(b)] $L^2 = \phi(L)^2 + \phi(L) - 2$ with $\phi(L) \geq 3$, $L \not\eqv 2D$ for $D$ such that 
$D^2 = 10$, $\phi(D) = 3$. In these cases $\gon(C) = 2\phi(L) - 1$ except for $\phi(L) = 3, 4$ when
$\gon(C) = 2\phi(L) - 2$.
\item[(c)] $(L^2, \phi(L)) = (30, 5)$, $(22, 4)$, $(20, 4)$, $(14, 3)$, $(12, 3)$ and $(6, 2)$. In 
these cases $\gon(C) = \lfloor \frac{L^2}{4} \rfloor + 2 = 2\phi(L) - 1$.
\end{itemize}
\end{bcor}

The line bundles in (a), (b) and (c) above have an explicit description by Proposition \ref{poscone}.

We also obtain the following result about the variation of the gonality of smooth curves in 
a complete linear system:

\begin{bcor}
\label{vargon}
Let $|L|$ be a base-component free complete linear system on an Enriques surface such that $L^2
> 0$. Let $\gengon |L|$ denote the gonality of a general smooth curve in $|L|$ and $\mingon |L|$ 
denote the minimal gonality of a smooth curve in $|L|$. Then
\[\gengon |L|-2 \leq \mingon |L| \leq \gengon |L|. \] 

Moreover if equality holds on the left, then $\phi(L) \geq \lceil \sqrt{\frac{L^2}{2}} \rceil$.
\end{bcor}

Finally in \S \ref{ex} we give examples showing that this result is sharp, that is that all the
cases $\mingon |L| = \gengon |L| - 2$, $\gengon |L| - 1$ and $\gengon |L|$ do occur.

Aside from our use of well-known vector bundle methods, to study linear series on
curves on an Enriques surface we will proceed, in section \ref{setting}, as follows. If a general
curve $C \in |L|$ carries some ``unexpected" linear series, then it also carries some zero-dimensional
schemes not imposing independent conditions on some subbundle of $L$. By moving these schemes on $S$ 
we will often derive a contradiction or find that the gonality is computed by $\mu(L)$. An important
feature of this method is that, unlike all the previous ones, it does work well, in many cases, also
{\it below the Bogomolov-Reider range}, that is when $L^2 < 4 \gon(C)$ (see Proposition
\ref{lemma:somma}), and also for other types of surfaces besides Enriques surfaces.

\vskip .2cm

\noindent {\it Acknowledgments}. We wish to thank Roberto Mu\~{n}oz for several helpful
discussions.

\section{Basic results on line bundles on Enriques surfaces}
\label{basic}

\begin{defn}
We denote by $\sim$ (resp. $\eqv$) the linear (resp. numerical) equivalence of divisors or line bundles. A
line bundle $L$ is {\bf primitive} if $L \eqv kL'$ implies $k = \pm 1$. We will write $L \geq 0$
(respectively $L > 0$) for an effective (resp. effective non trivial) line bundle $L$. If $V \subseteq H^0(L)$
is a linear system, we denote its {\bf base scheme} by $\Bs |V|$. We denote by $|L|_{sm}$ the open subset of
smooth curves in $|L|$. A {\bf nodal} curve on an Enriques surface $S$ is a smooth rational curve contained in
$S$.
\end{defn} 

We will often use another property of the function $\phi(L)$ defined in the introduction:
If $L$ is nef, then there exists a nef divisor $E$ calculating $\phi(L)$ (\cite[2.11]{co2} or by
\cite[Cor.2.7.1, Prop.2.7.1 and Thm.3.2.1]{cd}).

We also recall the following simple consequence of the signature theorem \cite[VIII.1]{bpv}.

\begin{lemma} \cite[Lemma 2.1]{klvan}
\label{lemma10}
Let $X$ be a smooth surface and let $A > 0$ and $B > 0$ be divisors on $X$ such that
$A^2 \geq 0$ and $B^2 \geq 0$. Then $A.B \geq 0$ with equality if and only if there exists a
primitive divisor $F > 0$ and integers $a \geq 1, b \geq 1$ such that $F^2 = 0$ and $A
\eqv aF, B \eqv bF$.
\end{lemma}

This gives

\begin{lemma} 
\label{A}
Let $L > 0$ and  $\Delta > 0$ be divisors on an Enriques surface $S$ with $L^2 \geq 0$, $\Delta^2 =
-2$ and $k := - \Delta.L > 0$. Then there exists an $A > 0$ such that $A^2 = L^2$, $A.\Delta=k$ and 
$L \sim A + k \Delta$. Moreover, if $L$ is primitive, then so is $A$. 
\end{lemma}

\begin{proof} Set $A = L - k \Delta$. Then $A.\Delta = k$ and $A^2 = L^2  \geq 0$. If $K_S - A \geq
0$, Lemma \ref{lemma10} yields $0 \leq (K_S-A).L = - L^2 - k^2$, a contradiction. Hence $H^2(A) = 0$ and
by Riemann-Roch we get $A > 0$. Now if $A \eqv qB$ with $q \geq 2$, then $k = qB.\Delta$, whence $L \eqv
q(B + (B.\Delta)\Delta)$ is not primitive.
\end{proof}

\begin{lemma} 
\label{lemma:nefred} 
Let $S$ be an Enriques surface and let $L$ be a line bundle on $S$ such that $L > 0, L^2 > 0$. Let
$F > 0$ be a divisor on $S$ such that $F^2 = 0$ and $\phi(L) = |F.L|$. Then $F.L > 0$ and if $\alpha > 0$ is
such that $(L - \alpha F)^2 \geq 0$, then $L - \alpha F > 0$.
\end{lemma}

\begin{proof}
By Lemma \ref{lemma10} we get $F.L > 0$. If $(L - \alpha F)^2 \geq 0$ we get by Riemann-Roch that either $L -
\alpha F > 0$ or $K_S - L + \alpha F > 0$. But in the latter case, Lemma \ref{lemma10} gives the contradiction
$- \phi(L) = F.(K_S - L + \alpha F) \geq 0$.
\end{proof}

We recall here a consequence of the vanishing theorem proved in \cite{klvan} that will be
used throughout the article.

\begin{defn} 
\label{def:qnef}
An effective line bundle $L$ on a K3 or Enriques surface is said to be {\bf quasi-nef} if $L^2
\geq 0$ and $L.\Delta \geq -1$ for every $\Delta$ such that $\Delta > 0$ and $\Delta^2 = -2$. 
\end{defn}

\begin{thm} 
\label{cor:qnef} \cite[Corollary 2.5]{klvan} 
An effective line bundle $L$ on a K3 or Enriques surface is quasi-nef if and only if $L^2 \geq 0$
and either $h^1(L) = 0$ or $L \eqv nE$ for some $n \geq 2$ and some primitive and nef divisor $E > 0$
with $E^2 = 0$.
\end{thm}
 
\subsection{Line bundles with $\mu(L) < 2 \phi(L)$} 

We want to prove a result about the function $\mu(L)$. In fact, the cases for 
which $\mu(L) < 2\phi(L)$ are of a very particular type and we will classify them.

\begin{defn}
Let $E_1, E_2, E_3$ be three  primitive divisors on an Enriques surface $S$ such 
that $E_i > 0$, $E_i^2 = 0$, $i = 1, 2, 3$, $E_1.E_2 = E_1.E_3 = 2$ and $E_2.E_3 = 1$.

A line bundle $L$ on $S$ is said to be of {\bf type}
\begin{itemize}
\item[{\bf ($\mu_1$)}] if $L \eqv h(E_1 + E_2)$, $h \geq 1$, 
\item[{\bf ($\mu_2$)}] if $L \sim h(E_1 + E_2) + E_3$, $h \geq 1$,
\item[{\bf ($\mu_3$)}] if $L \sim (h + 1)E_1 + hE_2 + E_3$, $h \geq 1$.
\end{itemize}
\end{defn} 

The properties of these three line bundle types will be proved in Lemma \ref{lemmaprop}.

Our goal here is to prove the ensuing

\begin{prop} 
\label{prop:mu} 
Let $L > 0$ be a line bundle on an Enriques surface $S$ such that $L^2 > 0$ and $(L^2, \phi(L))
\neq (4, 2)$. Then $\mu(L) < 2\phi(L)$ if and only if $L$ is of type ($\mu_1$), ($\mu_2$) or 
($\mu_3$).
\end{prop}

\begin{rem} The above proposition also holds when $(L^2, \phi(L)) = (4, 2)$. {\rm In fact in this
case we can prove that $\mu(L) = 3$. However the proof of this is quite long and will be omitted, as we
do not need it for the sequel.} 
\end{rem}
 
To prove the proposition we first need a few auxiliary results.

\begin{lemma} 
\label{lemma:lat2} 
Let $D > 0$ be a divisor on an Enriques surface $S$ such that $D^2 = 10$ and $\phi(D) = 3$. Then 
there are ten divisors $F_i$ such that $F_i > 0$, $F_i^2 = 0$ for $1 \leq i \leq 10$, $F_i.F_j =
1$ for $1 \leq i < j \leq 10$ and $3D \sim F_1 + \ldots + F_{10}$. Moreover if $F$ and $F'$
satisfy $F > 0$, $F' > 0$, $F^2 = (F')^2 = 0$, $F.D = F'.D = 3$ and $F \not\eqv F'$, then $F.F' =
1$.
\end{lemma}

\begin{proof}
The first assertion easily follows from \cite[Cor.2.5.5]{cd}, together with Lemma \ref{lemma10} for the
effectiveness of the $F_i$. For the last assertion, $F.(3D) = 9$ implies that $F.F_j = 0$ for some $j$, whence
$F \eqv F_{j_1}$, $F' \eqv F_{j_2}$ with $j_1 \neq j_2$ by Lemma \ref{lemma10}, so that $F.F' = 1$.
\end{proof}

\begin{lemma} 
\label{lemma:lat3}
Let $L > 0$ be a line bundle on an Enriques surface $S$ such that $L^2 > 0$. Then $L^2 = \phi(L)^2$ if and only
if $L$ is of type $\mu_1$.  
\end{lemma}

\begin{proof}
The ``if'' part is immediate. For the converse, assume $L^2 = \phi(L)^2$. Since $L^2$ is even,
$\phi(L)$ must be even, say $\phi(L) = 2h$, so that $L^2 = 4h^2$. The result is clear for $h = 1$, so we assume
henceforth that $h \geq 2$. Pick any $E_1 > 0$ with $E_1^2 = 0$ and $E_1.L = 2h$. Set $L_1 = L - E_1$, so that
$L_1^2 = 4h(h - 1)$ and $L_1 > 0$ by Lemma \ref{lemma:nefred}.

If $\phi(L_1) \geq 2h - 1$, then $4h^2 - 4h = L_1^2 \geq \phi(L_1)^2 \geq (2h - 1)^2$, a
contradiction. Therefore $\phi(L_1) \leq 2h - 2$. Pick any $E_2 > 0$ with $E_2^2 = 0$ and $E_2.L_1 =
\phi(L_1)$. Then $2h - E_1.E_2 \leq E_2.L - E_2.E_1 = E_2.L_1 \leq 2h - 2$, whence $E_1.E_2 \geq 2$ and if
equality holds we have that $E_2.L = 2h$. Now $(E_1 + E_2).L \leq 4h - 2 + E_1.E_2$, whence $(L_1 - E_2)^2 \geq
4(h - 1)^2 \geq 4$ and $L_1 - E_2 > 0$ by Lemma \ref{lemma:nefred}. The Hodge index theorem
on $E_1 + E_2$ and $L$ yields  that either $E_1.E_2 = 2$ or $E_1.E_2-2 \geq 8h(h - 1)$. In the latter case 
we have $E_1.(L_1 - E_2) \leq 2h - 8h(h - 1) - 2 < 0$, contradicting Lemma \ref{lemma10}.

Hence $E_1.E_2 = 2$, $E_1.L = E_2.L = 2h$ and we are done by the Hodge index theorem.
\end{proof}

Now we need to prove an integer version of \cite[Lemma 1.4.2]{co2}.

\begin{lemma} 
\label{lemma:lat1} 
Let $L > 0$ be a line bundle on an Enriques surface $S$ with $L^2 \geq 0$. Then there is an integer 
$n$ such that $1 \leq n \leq 10$ and, for every $i = 1, \ldots, n$, there are primitive divisors 
$E_i > 0$ with $E_i^2 = 0$ and integers $a_i > 0$ such that 
\[ L \eqv a_1E_1 + \ldots + a_nE_n \]
and one of the three following intersection sets occurs:
\begin{itemize}
\item[(i)] $E_i.E_j = 1$ for $1 \leq i < j \leq n$.
\item[(ii)] $n \geq 2$, $E_1.E_2 = 2$ and $E_i.E_j = 1$ for $2 \leq i < j \leq n$ and for $i = 1$,
$3 \leq j \leq n$.
\item[(iii)] $n \geq 3$, $E_1.E_2 = E_1.E_3 = 2$ and $E_i.E_j = 1$ for $3 \leq i < j \leq n$, for 
$i = 1$, $4 \leq j \leq n$ and for $i = 2$, $3 \leq j \leq n$.
\end{itemize} 
\end{lemma}

\begin{proof}
We prove the lemma by induction on $L^2$. Since the case $L^2 = 0$ is obvious, we assume $L^2 \geq 2$. By
Lemma \ref{lemma:nefred} we can choose a primitive divisor $F > 0$ with $F^2 = 0$, $F.L = \phi(L)$ and, setting
$L_1 = L - F$, we get that $0 \leq L_1^2 < L^2$ and $L_1 > 0$. By induction, we have that $L_1 \eqv a_1E_1 +
\ldots + a_nE_n$ and $L \eqv F + a_1E_1 + \ldots + a_nE_n$ with the intersections among the $E_i$'s as in (i),
(ii) or (iii). Note that if $F.E_i \leq  1$ for all $i = 1, \ldots, n$ then, by Lemma \ref{lemma10}, we have
that either $F \eqv E_i$ for some $i$ and then $L$ has the desired decomposition or $F.E_i = 1$ for all $i$. In
the latter case we cannot have $n = 10$ because the intersection matrix of $F, E_1, \ldots, E_{10}$ has
nonzero determinant while the Enriques lattice has rank $10$. Hence $n \leq 9$ and we are
done if $F.E_i \leq 1$ for all $i = 1, \ldots, n$. 

Therefore we will henceforth assume that $F.E_i \geq 1$ for all $i = 1, \ldots, n$ and 
that there is an index $i_0$ such that $F.E_{i_0} \geq 2$. 

We divide the proof in the three cases corresponding to the intersections of the 
$E_i$'s. To simplify some computations we set $a:= \sum\limits_{i=1}^n a_i$.

{\bf Case 1}: $E_i.E_j = 1$ for $1 \leq i < j \leq n$.

If $n \geq 2$ pick $j \in \{1, \ldots, n \} - \{i_0 \}$. Then $E_j.L = E_j.F + a
- a_j \geq \phi(L) = F.L \geq  a - a_j - a_{i_0} + a_{i_0}E_{i_0}.F + a_jE_j.F \geq a - a_j + 
a_{i_0} + a_jE_j.F$ giving the contradiction $0 \geq a_{i_0} + (a_j - 1)E_j.F \geq 1$. Therefore
$n = i_0 = 1$ and $E_1.F = E_1.L \geq \phi(L) = F.L = a_1E_1.F$ that is $a_1 = 1$. Now $\phi(L)^2 =
(E_1.F)^2 \leq L^2 = 2 E_1.F$, whence $E_1.F = 2$ and we are done for Case 1.

{\bf Case 2}: $n \geq 2$, $E_1.E_2 = 2$ and $E_i.E_j = 1$ for $2 \leq i < j \leq n$ and 
for $i = 1$, $3 \leq j \leq n$.

If $n \geq 3$ we have $E_3.L = E_3.F + a - a_3 \geq \phi(L) = F.L = 
\sum\limits_{i=1}^n a_i E_i.F$. If there is an $i_1 \in \{1, \ldots, n \} - \{3 \}$ such that
$F.E_{i_1} \geq 2$ then, as above, $F.L \geq a - a_3 + a_{i_1} + a_3E_3.F$, giving the
contradiction $0 \geq a_{i_1} + (a_3 - 1)E_3.F \geq 1$. Therefore we have that $F.E_i = 1$ for
all $i \in \{1, \ldots, n \} - \{3 \}$ and $i_0 = 3$. Then $E_3.L = E_3.F + a - a_3 \geq F.L \geq
a_3E_3.F + a - a_3$, whence $a_3 = 1$. Now $E_1.L = 1 + a - a_1 + a_2 \geq F.L \geq a + 1$ and
$E_2.L = 1 + a - a_2 + a_1 \geq F.L \geq a + 1$, giving $a_1 = a_2$ and $E_3.F = 2$. If $n \geq
4$ we have $E_4.L = 1 + a - a_4 \geq F.L = a + 1$, a contradiction. Hence $n = 3$ and $\phi(L) =
2a_1 + 2$, $L^2 = 4(a_1 + 1)^2$. Therefore $L^2 =\phi(L)^2$ and we are done by Lemma
\ref{lemma:lat3}.

If $n = 2$ we suppose, without loss of generality, that $E_2.F \geq 2$. From $E_i.L 
\geq F.L$, $i = 1, 2$, we get $0 \geq (a_1 - 1)E_1.F + a_2(E_2.F - 2)$ and $0 \geq (a_2 - 1)E_2.F
+ a_1(E_1.F - 2)$. The first inequality gives $a_1 = 1$, $E_2.F = 2$ and the second becomes $0
\geq 2(a_2 - 1) + (E_1.F - 2)$. If $E_1.F \geq 2$ we get that also $a_2 = 1$, $E_1.F = 2$, whence
$\phi(L) = F.L = 4$ while $L^2 = 12$, a contradiction. Therefore $E_1.F = 1$ and we are done in
this case.

{\bf Case 3}: $n \geq 3$, $E_1.E_2 = E_1.E_3 = 2$ and $E_i.E_j = 1$ for $3 \leq i 
< j \leq n$, for $i = 1$, $4 \leq j \leq n$ and for $i = 2$, $3 \leq j \leq n$.

If $n \geq 4$ we have $E_4.L = E_4.F + a - a_4 \geq \phi(L) = F.L = 
\sum\limits_{i=1}^n a_i E_i.F$. If there is an $i_1 \in \{1, \ldots, n \} - \{4 \}$ such that
$F.E_{i_1} \geq 2$ then, as above, $F.L \geq a - a_4 + a_{i_1} + a_4E_4.F$, giving the
contradiction $0 \geq a_{i_1} + (a_4 - 1)E_4.F \geq 1$. Therefore  we have that $F.E_i = 1$ for
all $i \in \{1, \ldots, n \} - \{4 \}$ and $i_0 = 4$. Then $E_4.(E_1 + E_2 + E_3) = F.(E_1 + E_2
+ E_3) = 3$, whence $E_4.F = 1$ by Lemma \ref{lemma:lat2}, a contradiction. 

Hence $n = 3$ and from $E_i.L \geq F.L$, $i = 1, 2, 3$, we get
\begin{eqnarray}
\label{eq1} & 0 \geq (a_1 - 1)E_1.F + a_2(E_2.F - 2) + a_3(E_3.F - 2), \\ 
\label{eq2} & 0 \geq a_1(E_1.F - 2) + (a_2 - 1)E_2.F + a_3(E_3.F - 1), \\ 
\label{eq3} & 0 \geq a_1(E_1.F - 2) + a_2(E_2.F - 1) + (a_3 - 1)E_3.F.
\end{eqnarray}
Now if $E_1.F \geq 2$ we get from (\ref{eq3}) that $E_1.F = 2, E_2.F = a_3 = 1$. Then (\ref{eq2}) 
gives $E_3.F = a_2 = 1$ and (\ref{eq1}) implies that $1 \leq a_1 \leq 2$. Then $\phi(L) = F.L =
2a_1 + 2$ while $L^2 = 12a_1 + 6$. This gives a contradiction when $a_1 = 2$ since then $\phi(L) = 6$
and $L^2 = 30 < 36$. Therefore $a_1 = 1$, $\phi(L) = 4$ and $L^2 = 18$. Now $(L - 2F)^2 = 2$
whence, by Lemma \ref{lemma:nefred}, we can write $L - 2F \sim F_1 + F_2$ with $F_i > 0$, $F_i^2 =
0$, $i = 1, 2$ and $F_1.F_2 = 1$. Also $4 \leq F_i.L = 2F_i.F + 1$, therefore $F_i.F \geq 2$ for
$i = 1, 2$. Since $F_1.F + F_2.F = F.(L -2F) = 4$ we get $F_1.F = F_2.F = 2$ and we are done in
this case.

Therefore we can assume in the sequel of the proof that $E_1.F = 1$. 

Now if $E_i.F \geq 2$ for $i = 2, 3$ we get from (\ref{eq1}) that $E_2.F = E_3.F = 2$ and
$a_1 = 1$. Adding up (\ref{eq2}) and (\ref{eq3}) gives $0 \geq -6 + 3a_2 + 3a_3 \geq 0$, therefore
also $a_2 = a_3 = 1$. But then $\phi(L) = F.L = 5$ and $L^2 = 20$, a contradiction.

Therefore we can assume, without loss of generality, that $E_1.F = E_2.F = 1$, $E_3.F 
\geq 2$. 

Now (\ref{eq1}) becomes $a_2 \geq a_1 + a_3(E_3.F - 2) - 1$ and (\ref{eq2}) becomes 
$a_1 \geq a_2 + a_3(E_3.F - 1) - 1$. Combining we get $a_3 (2E_3.F - 3) \leq 2$, whence $E_3.F = 2$
and $1 \leq a_3 \leq 2$. Using again the inequalities (\ref{eq1}) and (\ref{eq2}) we get the only
possibilities $a_3 = 2, a_2 = a_1 - 1$ or $a_3 = 1, a_1 - 1 \leq a_2 \leq a_1$. In the first case we get
$\phi(L)^2 = (F.L)^2 = (2a_1 + 3)^2 > L^2 = 4a_1^2 + 12a_1 + 2$, a contradiction. In the second case,
setting $b = a_1$, we have the two possibilities
\begin{eqnarray}
\label{eqb} & L \eqv F + b E_1 + b E_2 + E_3 \\ 
\label{eqb-1} & L \eqv F + b E_1 + (b - 1) E_2 + E_3.
\end{eqnarray}
Set $D = E_1 + E_2 + E_3$ so that $D^2 = 10$ and $\phi(D) = 3$. We can write $3D \sim F_1 + \ldots +
F_{10}$ as in Lemma \ref{lemma:lat2}. Since $E_2.D = E_3.D = 3$ by Lemma \ref{lemma10} we can assume, without
loss of generality, that $E_2 \eqv F_1, E_3 \eqv F_2$, whence $E_2.F_i = E_3.F_i = 1$ for $3 \leq i \leq 10$.
Also $F.D = 4$, therefore, by Lemma \ref{lemma10}, $F.F_i \geq 1$ for $1 \leq i \leq 10$, and we can also
assume, without loss of generality, that $F.F_3 = 2$. Also from $12 = 3E_1.D = 4 + E_1.F_3 + \ldots +
E_1.F_{10}$ we see that $E_1.F_3 = 1$. Now set $D' = F + F_3 + E_3$ so that $(D')^2 = 10$ and $D.D' = 10$,
therefore $D \eqv D'$ by the Hodge index theorem. It follows that $E_1 + E_2 \eqv F + F_3$.

If $L$ is as in (\ref{eqb}) then $L \eqv (b + 1) F + b F_3 + E_3$ with $F.F_3 = F.E_3 
= 2$, $F_3.E_3 = 1$, has the required decomposition. 

Now suppose that $L$ is as in (\ref{eqb-1}). Set $F' = L - bF - b F_3$. Now $(F')^2 = 0,
F.F' = 1$ whence $F' > 0$ by Riemann-Roch. Also $F_3.F' = 2$ therefore $L \sim b F + b F_3 + F'$ with
$F_3.F = F_3.F' = 2$, $F.F' = 1$, has the required decomposition. 
\end{proof}

\begin{lemma} 
\label{lemma:lat4} 
Let $L > 0$ be a line bundle on an Enriques surface $S$ with $L^2 > 0$. Then $\mu(L) \geq 2\phi(L)
-  2$. Moreover if $\mu(L) < 2\phi(L)$ and $B$ is a line bundle computing $\mu(L)$, that is $B > 0$,
$B^2 = 4$, $\phi(B) = 2$, $B \not\eqv L$ and $B.L = \mu(L)$, then $B \sim F_1 + F_2$, with $F_i >
0$ primitive, $F_i^2 = 0$, $i = 1, 2$, $F_1.F_2 = 2$, $F_1.L = \phi(L)$ and $F_2.L =
\phi(L)$ or $\phi(L) + 1$.
\end{lemma}

\begin{proof}
Use Lemma \ref{lemma:nefred}.
\end{proof}

\begin{lemma} 
\label{lemmaprop}
\null \hskip 7cm
\begin{itemize}
\item[(i)] If $L$ is of type ($\mu_1$) then $\phi(L) = 2h$ and $L^2 = \phi(L)^2$. If $h \geq 2$ then
$\mu(L) = 2\phi(L) - 2$.
\item[(ii)] If $L$ is of type ($\mu_2$) then $\phi(L) = 2h + 1$, $\mu(L) = 2\phi(L) - 1$ and $L^2
= \phi(L)^2 + \phi(L) - 2$.
\item[(iii)] If $L$ is of type ($\mu_3$) then $\phi(L) = 2h + 2$, $\mu(L) = 2\phi(L) - 1$ and $L^2
= \phi(L)^2 + \phi(L) - 2$.
\end{itemize}
\end{lemma}

\begin{proof}
Apply Lemmas \ref{lemma:lat2}, \ref{lemma:lat3} and \ref{lemma:lat4}. 
\end{proof}

\renewcommand{\proofname}{Proof of Proposition {\rm \ref{prop:mu}}} 
\begin{proof}
By Lemma \ref{lemmaprop} we can assume that $L^2 > 0$, $(L^2, \phi(L)) \neq (4, 2)$ and $\mu(L) <
2\phi(L)$. We want to  show that $L$ must be of type ($\mu_1$), ($\mu_2$) or ($\mu_3$). To do this, we divide
the treatment into the three cases occurring in Lemma \ref{lemma:lat1}. 

We set $a = \sum\limits_{i = 1}^n a_i$ and we choose $B \sim F_1 + F_2$ as in Lemma \ref{lemma:lat4}
that computes $\mu(L)$, so that $F_1.L = \phi(L)$ and $F_2.L = \phi(L)$ or $\phi(L) + 1$.

If $L$ is as in (i) of Lemma \ref{lemma:lat1}, then reordering the $a_i$'s so that 
$a_1 \geq \ldots \geq a_n$ we have $a - a_1 = E_1.L \leq \ldots \leq E_n.L = a - a_n$. Now any $F
> 0$ with $F^2 = 0$, $F \not\eqv E_i$ for all i, satisfies $F.L \geq a$. Hence
$\phi(L) = a - a_1 < a$. Now $F_1.L = a - a_1$, whence we can assume $F_1 \eqv E_1$ after renumbering
indices. Since $F_1.F_2 = 2$ we have $F_2 \not\eqv E_i$ for $i \geq 2$, whence $F_2.L \geq a + a_1$, so that
$\mu(L) = B.L - 2 \geq 2(a - 1) \geq 2 \phi(L)$, a contradiction. 

If $L$ is as in (ii) of Lemma \ref{lemma:lat1} we are done if $n = 2$ and $a_1 = a_2$. We
assume that this is not the case. Reordering the $a_i$'s so that $a_1 \geq a_2$ and $a_3 \geq \ldots
\geq a_n$ we have $a + a_2 - a_1 = E_1.L \leq E_2.L = a +  a_1 - a_2$ and $a - a_3 = E_3.L \leq \ldots
\leq E_n.L = a - a_n$. Now any $F > 0$ with $F^2 = 0$, $F \not\eqv E_i$ for all $i$, satisfies $F.L \geq a$.
Hence $\phi(L) = E_1.L = a + a_2 - a_1 \leq a$ if $n = 2$ and $\phi(L) = \min \{E_1.L, E_3.L \} = \min \{a +
a_2 - a_1, a - a_3 \} < a$ if $n \geq 3$. Since $F_1.L = \phi(L)$, we can assume, after renumbering indices,
that $F_1 \eqv E_1$ or $F_1 \eqv E_3$. 

If $F_1 \eqv E_1$, then $\phi(L) = a + a_2 - a_1$ and either $F_2 \eqv E_2$ or $F_2 \not\eqv
E_i$ for all $i$. In the first case, by Lemma \ref{lemma:lat4}, we have $2(a + a_2 -
a_1) + 1 = 2\phi(L) + 1 \geq (F_1 + F_2).L = (E_1 + E_2).L = 2a$, whence $a_1 = a_2$ and $n = 2$, a
contradiction. In the second case we have $F_2.L \geq a + a_1$, whence $2(a + a_2 - a_1) + 1 = 2\phi(L)
+ 1 \geq (F_1 + F_2).L = (E_1 + F_2).L \geq 2a + a_2$, therefore $a_2 + 1 \geq 2a_1 \geq 2a_2$, that
gives $a_1 = a_2 = 1$. Then $\phi(L) = a$, whence $n = 2$, again excluded.

If $F_1 \eqv E_3$, then $\phi(L) = a - a_3 < a$ and $F_2 \not\eqv E_i$ for all $i$. Therefore $F_2.L
\geq a + a_3 \geq \phi(L) + 2$, contradicting Lemma \ref{lemma:lat4}.

Finally, if $L$ is as in (iii) of Lemma \ref{lemma:lat1}, we claim that we can write
\begin{eqnarray}
\label{eq:lat1}
L \eqv a_1E_1 + \ldots + a_{10}E_{10}, \; \mbox{with} \; a_1 > 0, \; a_2 \geq a_3 > 0, \; 
a_{4} \geq \ldots \geq a_{10} \geq 0 \\
\nonumber  \; E_1.E_2 = E_1.E_3 = 2
\; \mbox{and} \; E_i.E_j = 1 \; \mbox{for} \; (i,j) \neq (1,2),(2,1),(1,3),(3,1).
\end{eqnarray}
Indeed, if $n = 10$, this is clear after renumbering indices. If $3 \leq n < 10$, we note that
$D := E_1 + E_2 + E_3$ satisfies $D^2 = 10$ and $\phi(D) = 3$, so that $3D \sim F_1 + \ldots + F_{10}$ as in
Lemma \ref{lemma:lat2}. Since $E_i.(3D) = 9$ for all $i = 2, \ldots, n$ by assumption, we must have $E_i.F_j =
0$ for some $j \in \{1, \ldots, 10\}$ for any $i = 2, \ldots, n$, whence for $2 \leq i \leq n$ we have
$E_i$ or $E_i + K_S \in \{F_1, \ldots, F_{10} \}$, so that each $F_i \not\eqv E_j$ for $j = 2 ,3$, satisfies
$F_i.E_1 = 1$. Therefore we can complete the set $\{ E_1,
\ldots, E_n \}$ to a set $\{E_1, \ldots, E_n, E_{n+1}, \ldots, E_{10} \}$ satisfying the
desired conditions (setting $a_{n+1} = \ldots = a_{10} = 0$).

Now, using (\ref{eq:lat1}), we deduce $a - a_4 = E_4.L \leq \ldots \leq E_{10}.L = a - a_{10}$,
$a - a_2 + a_1 = E_2.L \leq a - a_3 + a_1 = E_3.L$ and $E_1.L = a + a_2 + a_3 - a_1$. Moreover, for any $F > 0$
with $F^2 = 0$ and $F \not\eqv E_i$, $1 \leq i \leq 10$, we have $F.L \geq a$. 
Combining the above we see that
\begin{equation} 
\label{eq:lat5}
\phi(L) = \min \{ E_1.L, E_2.L, E_4.L \} = \min \{ a + a_2 + a_3 - a_1, a - a_2 + a_1, a - a_4 \}.
\end{equation} 
Next we will prove that $\mu(L) \geq 2a + a_3 - 2$. Suppose then that $\mu(L) < 2a + a_3 - 2$.

Assume first that $F_i \not\eqv E_j$ for $i = 1, 2$ and all $j$. Then $F_i.L \geq a$ and using
\eqref{eq:lat5} we have 
\[ 2a - 2 \leq \mu(L) < 2\phi(L) = \min \{ 2(a + a_2 + a_3 - a_1), 2(a - a_2 + a_1), 2(a - a_4) 
\}. \]
It follows that $a_4 = 0$ (whence $a_5 = \ldots = a_{10} = 0$) and
\begin{equation} 
\label{eq:lat8}
a_2 \leq a_1 \leq a_2 + a_3.
\end{equation} 
Moreover $\phi(L) = a$ and $\mu(L) = 2a - 2$ or $2a - 1$, so that $F_1.L = a$ and $F_2.L = a$ 
or $ a + 1$. In the first case we have $F_1.(E_1 + E_2 + E_3) = F_2.(E_1 + E_2 + E_3) = 3$, and
from Lemma \ref{lemma:lat2} we get $F_1.F_2 = 1$, a contradiction. In the second case, using our
assumption that $\mu(L) < 2a + a_3 - 2$, we must have $a_3 \geq 2$, and consequently $a_2 \geq 2$. Hence
$a_1 \geq 2$ as well by \eqref{eq:lat8}. Therefore $F_2.L = a_1F_2.E_1 + a_2F_2.E_2 + a_3F_2.E_3 = a$ or
$\geq a + 2$, a contradiction.

Assume next that $F_1 \eqv E_i$ for some $i$ but $F_2$ is not. Then, using \eqref{eq:lat5}, we can
assume $F_1 \eqv E_1$, $E_2$ or $E_4$, after renumbering indices (but still maintaining the inequalities in
(\ref{eq:lat1})).

If $F_1 \eqv E_1$ then $F_1.L = \phi(L) = a + a_2 + a_3 - a_1$ and $F_2.L \geq a + a_1$,
whence $\mu(L) \geq 2a + a_2 + a_3 - 2 \geq 2a + 2a_3 - 2 > 2a + a_3 - 2$, a contradiction.

If $F_1 \eqv E_2$ then $F_1.L = \phi(L) = a - a_2 + a_1$ and $F_2.L \geq a + a_2$, 
whence $\mu(L) \geq 2a + a_1 - 2$. From \eqref{eq:lat5} we get
\[ 2a + a_1 - 2 \leq \mu(L) < 2\phi(L) =  2(a - a_2 + a_1), \]
whence $2a_2 \leq a_1 + 1$ and $a - a_2 + a_1 \leq a - a_4$, so that $a_4 \leq a_2 - a_1 \leq 
1 - a_2 \leq 0$. Therefore $a_1 = a_2 = a_3 = 1$, $a_4 = a_5 = \ldots = a_{10} = 0$ and we get
the contradiction $5 = 2a + a_1 - 2 \leq \mu(L) < 2a + a_3 - 2 = 5$.

If $F_1 \eqv E_4$ then $F_1.L = \phi(L) = a - a_4$ and $a - a_4 + 1 = \phi(L) + 1 
\geq F_2.L \geq a + a_4$, so that $a_4 = a_5 = \ldots = a_{10} = 0$ and $\phi(L) = a$, $\mu(L) =
2a - 2$ or $2a - 1$. Moreover, by \eqref{eq:lat5} we see that \eqref{eq:lat8} holds and we
derive the same contradiction as above, right after \eqref{eq:lat8}. 

Now assume that $F_2 \eqv E_i$ for some $i$ but $F_1$ is not. Then $a \leq F_1.L = \phi(L) \leq a -
a_4$ whence $a_4 = \ldots = a_{10} = 0$ and $F_1.E_i = 1$ for $1 \leq i \leq 3$. Since $F_1.F_2 = 2$ we get $F_2
\eqv E_j$ for some $j \geq 4$, but then $F_1.(E_1 + E_2 + E_3) = F_2.(E_1 + E_2 + E_3) = 3$, contradicting
Lemma \ref{lemma:lat2}. 

We have therefore proved that if $\mu(L) < 2a + a_3 - 2$ then $F_1 \eqv E_i$, $F_2 \eqv E_j$ for some
$i, j$. By \eqref{eq:lat1} and \eqref{eq:lat5} we can assume $F_1 \eqv E_1$ and $F_2 \eqv E_2, E_3$ or $F_1
\eqv E_2$ and $F_2 \eqv E_1$. Now $(E_1 + E_2).L = 2a + a_3 \leq 2a + a_2 = (E_1 + E_3).L$, contradicting
$\mu(L) < 2a + a_3 - 2$. 

This proves that $\mu(L) \geq 2a + a_3 - 2$. Comparing with \eqref{eq:lat5} we have
\[2a + a_3 - 2 \leq \mu(L) < 2\phi(L) = \min \{ 2(a + a_2 + a_3 - a_1), 2(a - a_2 + a_1), 
2(a - a_4) \}. \]
Then $2a + a_3 - 2 < 2a - 2 a_4$, whence $a_3 = 1$, $a_4 = \ldots = a_{10} = 0$ and $a_2 \leq a_1
\leq a_2 + 1$, so that $L$ is of type ($\mu_2$) or ($\mu_3$). 
\end{proof}
\renewcommand{\proofname}{Proof}

\subsection{A result on the positive cone of an Enriques surface}

\renewcommand{\proofname}{Proof of Proposition {\rm \ref{poscone}}}  
\begin{proof}
We prove the proposition by induction on $L^2$. It is easily seen to hold fo $\phi(L) \leq 3$. If
$\phi(L) = 4$ we have that either $L^2 = 16$ and we get case (i) by Lemma \ref{lemma:lat3} or $L^2 = 18$.
In the latter case let $E > 0$ be such that $E^2 = 0$ and $E.L = 4$. 
Write $L - 2E \sim E_1 + E_2$ for $E_i > 0$ with $E_i^2 = 0$, $i = 1, 2$ and $E_1.E_2 = 1$ by Lemma
\ref{lemma:nefred}. Since $4 = \phi(L) \leq E_i.L = 2E.E_i + 1$, we find $E.E_i = 2$ for $i = 1,2$, whence $L$
is as in case (ii-b).

We will henceforth assume that $\phi(L) \geq 5$ and therefore $L^2 \geq 26$.

Pick an $E > 0$ such that $E^2 = 0$ and $E.L = \phi(L)$. Then $(L - E)^2 > 0$ and 
$L - E > 0$ by Lemma \ref{lemma:nefred}. Moreover $\phi(L) = E.L = E.(L - E) \geq \phi(L - E)$.
If $\phi(L - E) = \phi(L)$, then 
\[ (L - E)^2 = L^2 - 2\phi(L) \leq \phi(L)^2 - \phi(L) - 2 < \phi(L)^2 = \phi(L - E)^2, \]
a contradiction. 
If $L^2 = \phi(L)^2$ and $\phi(L - E) = \phi(L) - 1$, then $(L - E)^2 =  \phi(L - E)^2 - 1$, 
again a contradiction. We have therefore proved that 
\begin{equation} 
\label{eq:bassi1}
\phi(L - E) \leq \phi(L) -1 \; \mbox{with equality only if} \; L^2 > \phi(L)^2.
\end{equation}

Now pick an $E' > 0$ such that $(E')^2 = 0$ and $E'.(L - E) = \phi(L - E)$.

Assume $\phi(L - E) \leq \phi(L) - 2$. From $\phi(L) \leq E'.L = \phi(L - E) + E.E'$, we find $E.E'
\geq 2$. If equality holds, then Proposition \ref{prop:mu} and Lemma \ref{lemmaprop} imply that $L$ is as in
case (i). 

If $3 \leq E.E' \leq 6$, then $6 \leq (E + E')^2 \leq 12$ and it is easily shown that $E +
E'$ is  a sum of at least three divisors $F_i > 0$ with $F_i^2 = 0$, whence $(E + E').L \geq 3\phi(L)$.
Now $E'.L \leq \phi(L) - 2 + E.E' \leq \phi(L) + 4$, whence $\phi(L) \leq 4$, a contradiction.

If $E.E' \geq 7$ we set $\alpha = \lfloor \frac{\phi(L)}{2} \rfloor - 1$. Then
$(L - E - \alpha E')^2 \geq  2\phi(L) - 4 \geq 6$,
so that $L - E - \alpha E' > 0$ by Lemma \ref{lemma:nefred}, whence
\[ \phi(L) = E.L = E.(L - E - \alpha E') + \alpha E.E' \geq 1 + 7\alpha \geq 1 + 
7(\frac{\phi(L) - 3}{2}), \]
giving the contradiction $\phi(L) \leq 3$.

We have therefore proved that one of the following holds:
\begin{eqnarray} 
\label{eq:bassi2}
\hskip .8cm \phi(L - E) = \phi(L) - 1 \; \mbox{and} \; L^2 > \phi(L)^2,
\end{eqnarray}
\begin{eqnarray} 
\label{eq:bassi2'} 
& \phi(L - E) = \phi(L) - 2  \; \mbox{and} \;  L^2 = \phi(L)^2. 
\end{eqnarray}

Of course in case (\ref{eq:bassi2'}) we get that $L$ is as in (i) by
Lemma \ref{lemma:lat3}.
 
Now assume that we are in case (\ref{eq:bassi2}). Then 
\[ (L - E)^2 = L^2 - 2\phi(L) \leq \phi(L)^2 - \phi(L) - 2 = \phi(L - E)^2 + \phi(L - E) - 2. \]
If $(L - E)^2 = \phi(L - E)^2$, then, by Lemma \ref{lemma:lat3}, we have that $L - E \eqv 
h(E_1 + E_2)$ with $h \geq 1$, $E_i > 0$, $E_i^2 = 0$ and $E_1.E_2 = 2$. In particular 
$\phi(L - E) = 2h$, so that $\phi(L) = 2h + 1$. It follows that $h \geq 2$. Now $\phi(L) = 2h + 1 =
E.L = h(E.E_1 + E.E_2)$, whence we must have either $E.E_1 \geq 2$ or $E.E_2 \geq 2$. Moreover
$E.E_i > 0$ for $i = 1, 2$ by Lemma \ref{lemma10} as $E.(L - E) > E_i.(L - E)$. Hence
$\phi(L) \geq 3h$, a contradiction. Therefore we have
\begin{equation} 
\label{eq:bassi3}
\phi(L - E)^2 <  (L - E)^2 \leq \phi(L - E)^2 + \phi(L - E) - 2,
\end{equation}
and we can assume by induction that in fact $(L - E)^2 = \phi(L - E)^2 + \phi(L - E) - 2$ (and,
consequently, $L^2 = \phi(L)^2 + \phi(L) - 2$) and that we are in one of the three following cases,
where $h \geq 1$, all the $E_i$'s are primitive, $E_i > 0$, $E_i^2 = 0$,
$E_1.E_2 = E_1.E_3 = 2$ and $E_2.E_3 = 1$:
\begin{itemize}
\item[(a)] $L - E \sim h(E_1 + E_2) + E_3$, $\phi(L - E) = 2h + 1$,
\item[(b)] $L - E \sim (h + 1)E_1 + hE_2 + E_3$, $\phi(L - E) = 2h + 2$,
\item[(c)] $L - E \eqv 2(E_1 + E_2 + E_3)$, $\phi(L - E) = 6$.
\end{itemize}

Now case (c) cannot occur since we have $E.L = \phi(L) = \phi(L - E) + 1 = 7$.

In case (a) we have $E.(L - E) = \phi(L) = \phi(L - E) + 1 = 2h + 2$, whence $h \geq
2$, $E_1.(L - E) = 2h + 2$, $E_2.(L - E) = 2h + 1$ and $E_3.(L - E) = 3h$. Since $E$ and all the
$E_i$'s are primitive, we must have $E.E_2 > 0$ by Lemma \ref{lemma10}. For the same
reason, if $E.E_3 = 0$ we must have $E \eqv E_3$ and $h = 2$. Then $L \eqv 2(E_1 + E_2 + E_3)$ and
we are in case (ii-c). Again, if $E.E_1 = 0$ we must have $E \eqv E_1$, whence $L \eqv (h + 1)E_1 + hE_2 + E_3$
and we are in case (ii-b).

Therefore we can assume $E.E_i > 0$ for all $i = 1, 2, 3$.

If $E.E_1 \geq 2$ or $E.E_2 \geq 2$, then $2h + 2 = E.L \geq 3h + 1$, a contradiction. 
Hence $E.E_1 = E.E_2 = 1$ and $E.E_3 = 2$. Then $E_1.L = 2h + 3$, $E_2.L = 2h + 2$ and
$\mu(L) \leq (E_1 + E_2).L - 2 = 4h + 3 < 2\phi(L)$ and we are in case (ii-b) by Proposition \ref{prop:mu}
and Lemma \ref{lemmaprop}.

In case (b), working as in case (a), we deduce that $E.E_1 > 0$ and either $E \eqv E_3$, $h = 1$ and
$L \eqv 2E_1 + 2E_3 + E_2$ or $E \eqv E_2$ and $L \eqv (h + 1)(E_1 + E_2) + E_3$ and we are in case  (ii-a).
Therefore we can assume $E.E_i > 0$ for all $i = 1, 2, 3$.

If $E.E_1 \geq 2$, then $2h + 3 = E.L \geq 3h + 3$, a contradiction. 
If $E.E_2 \geq 2$, then $2h + 3 = E.L \geq 3h + 2$, so that we must have $h = 1$,
$E.E_2 = 2$ and $E.E_1 = E.E_3 = 1$. But then $E_1.L = 5$ and $E_3.L = 6$, so that $\mu(L) \leq 9 
< 2\phi(L) = 10$, and we must be in case (ii-a) by Proposition \ref{prop:mu} and
Lemma \ref{lemmaprop}.

Hence $E.E_1 = E.E_2 = 1$, $E_1.L = 2h + 3$ and $E_2.L = 2h + 4$, whence
$\mu(L) \leq (E_1 + E_2).L - 2 = 4h + 5 < 2\phi(L)$ and we are in case (ii-a) by Proposition
\ref{prop:mu} and Lemma \ref{lemmaprop}.
\end{proof}
\renewcommand{\proofname}{Proof}

\begin{rem} {\rm Proposition \ref{poscone} improves the result of Hana
\cite[Thm.1.8]{ha}}. 
\end{rem}

\begin{lemma} 
\label{lemma:bassi2}
If $L \eqv 2D$ with $D^2 = 10$ and $\phi(D) = 3$, then $\mu(L) = 2\phi(L) = 12$.
\end{lemma}

\begin{proof}
Apply Lemmas \ref{lemma10}, \ref{lemma:nefred} and \ref{lemmaprop}
and Proposition \ref{prop:mu}.
\end{proof}

\subsection{A few useful applications} 

A direct application of Lemmas \ref{lemma10}-\ref{lemma:nefred} and Theorem \ref{cor:qnef} 
yields the results in Lemmas \ref{cl:g=6,2}-\ref{lemma:g=9,3}, which will be of use to us.

\begin{lemma} 
\label{cl:g=6,2}
Let $|L|$ be a base-point free complete linear system on an Enriques surface $S$ with $L^2 = 10$ and
$\phi(L) = 2$. Let $E$ be a nef divisor such that $E^2 = 0$ and $E.L = 2$. Then
\[ L \sim 2E + E_1 + E_2 \]
for $E_i > 0$ primitive with $E_i^2 = 0$ and $E.E_i = E_1.E_2 = 1$, $i = 1, 2$ such that
\begin{itemize}
\item[(i)] $|E + E_1|$ is base-component free with two distinct base points,
\item[(ii)] $E + E_2$ is quasi-nef.
\end{itemize}
\end{lemma}

\begin{proof}
The existence of the  decomposition of $L$ follows from Lemmas \ref{lemma:nefred} and
\ref{lemma10}.

Assume there is a $\Delta > 0$ such that $\Delta^2 = -2$ and $\Delta.(E + E_1) < 0$. Then
$\Delta.E_1 < 0$ by the nefness of $E$. By Lemma \ref{A} there is an $A > 0$ primitive such that 
$A^2 = 0$ and $E_1 \sim A + k \Delta$, where $k = - \Delta.E_1 \geq 1$. Since
$1 = E.E_1 = E.A + k E.\Delta \geq k E.\Delta$ we get that if $\Delta.E > 0$ then $k = \Delta.E =
1$, whence $\Delta.(E + E_1) = 0$, a contradiction. Therefore $\Delta.E = 0$, whence $\Delta.E_2 \geq
k$ by the nefness of $L$ and since $1 = E_2.E_1 = E_2.A + k E_2.\Delta$ we get that 
$k = E_2.\Delta = 1$ and $A.E_2 = 0$, so that $A \eqv E_2$ by Lemma \ref{lemma10}. Hence $E_1 \eqv E_2 +
\Delta$, with $E_2.\Delta = 1$.  In particular $E + E_1$ is quasi-nef. Also, if $E + E_1$ is not nef
then there is a nodal curve $\Gamma$ such that $E_1 \eqv E_2 + \Gamma$, with $E_2.\Gamma = 1$.

Similarly $E + E_2$ is quasi-nef and if there is a nodal curve $\Gamma'$ such that 
$\Gamma'.(E + E_2) < 0$ then $E_2 \eqv E_1 +\Gamma'$, with $E_1.\Gamma' = 1$. 

Obviously it follows that either $E + E_1$ or $E + E_2$ is nef. By symmetry, we can assume that $E +
E_1$ is nef, and  by \cite[Prop.3.1.6, Cor.3.1.4 and Thm.4.4.1]{cd} the lemma is proved, 
possibly after adding $K_S$ to both $E_1$ and $E_2$. 
\end{proof}

\begin{lemma} 
\label{cl:g=6,3}
Let $|L|$ be a base-point free complete linear system on an Enriques surface $S$ with $L^2 = 10$ and
$\phi(L) = 3$. Among all $E > 0$ satisfying $E^2 = 0$ and $E.L = 3$ pick one which is maximal. Then
\[ L \sim E + E_1 + E_2 \]
for $E_i > 0$ primitive with $E_i^2 = 0$, $i = 1, 2$, $E.E_1 = 1$, $E.E_2 = E_1.E_2 = 2$ 
such that
\begin{itemize}
\item[(i)] $|E_1 + E_2|$ and $|E_1 + E_2 + K_S|$ are base-point free,
\item[(ii)] $h^1(2E - L) \leq 1$.
\end{itemize}
\end{lemma}

\begin{proof}
The existence of the  decomposition of $L$ follows from Lemmas \ref{lemma:nefred} and \ref{lemma10}. If $E_1
+ E_2$ is not nef, it follows that there is a nodal curve $\Gamma$ such that $(E + \Gamma)^2 = 0$ and $(E
+ \Gamma).L = 3$, contradicting the maximality of $E$. This proves (i). Note that $h^0(2E - L) = 0$
and $h^0(L - 2E + K_S) \leq 1$ since $L.(L - 2E + K_S) = 4 < 2\phi(L)$, giving (ii) by Riemann-Roch. 
\end{proof}

\begin{lemma} 
\label{lemma:g=8,2}
Let $|L|$ be a base-point free complete linear system on an Enriques surface $S$ with $L^2 = 14$
and $\phi(L) = 2$ and let $E > 0$ be a nef divisor with $E^2 = 0$ and $E.L = 2$. 
Then there exists a decomposition
\[ L \sim 3E + E_1 + E_2 \]
with $E_i > 0$, $E_i^2 = 0$ and $E.E_i = E_1.E_2 = 1$, $i = 1, 2$. Moreover $3E + E_1$ is nef and 
$h^0(E_2 + K_S) = 1$.
\end{lemma}

\begin{proof}
We only prove that $h^0(E_2+K_S)=1$. If $h^0(E_2 + K_S) \geq 2$ then $h^1(E_2 + K_S) \geq 1$ by
Riemann-Roch, whence, using  Theorem \ref{cor:qnef}, there exists a divisor $\Delta > 0$ such that $\Delta.E_2
\leq - 2$. By Lemma \ref{A} we can write $E_2 + K_S \sim A + k \Delta$ with $A > 0$, $A^2 = 0$ and $k = -
E_2.\Delta \geq 2$. If $E.\Delta > 0$ we get the contradiction $1 = E.E_2  \geq 2$. Therefore $E.\Delta = 0$
and  similarly $E_1.\Delta \leq 0$, contradicting the nefness of $L$.
\end{proof}

\begin{lemma} 
\label{lemma:g=9,3} 
Let $|L|$ be a base-point free complete linear system on an Enriques surface $S$ with $L^2 = 16$ and
$\phi(L) = 3$ and let $E > 0$ be a nef divisor with $E^2 = 0$ and $E.L = 3$. Then there exists a
decomposition
\[ L \sim 2E +E_1+E_2 \]
for $E_i > 0$ with $E_i^2 = 0$, $i = 1, 2$, $E.E_1 = E_1.E_2 = 2$, $E.E_2 = 1$ and
$2E + E_1$ is nef. Moreover, either  $h^0(2E + E_1 -  E_2) = h^0(2E + E_1 - E_2 + K_S) = 1$ or 
$E_1 \eqv 2F$, for an $F > 0$ with $F^2 = 0$.
\end{lemma}

\begin{proof}
The nefness of $2E + E_1$ follows, as in the previous lemmas, by choosing a maximal $E_2$. Now note that
$(2E + E_1 - E_2)^2 = 0$ and $E.(2E + E_1 - E_2) = 1$, whence $2E + E_1 - E_2 > 0$. If $h^0(2E + E_1 - E_2) > 1$
or if $h^0(2E + E_1 - E_2 + K_S) > 1$, then by Riemann-Roch and Theorem \ref{cor:qnef}, there must exist a
$\Delta > 0$ with $\Delta^2 = -2$ and $\Delta.(2E + E_1 - E_2) \leq -2$. Since $2E + E_1$ is nef, we  must have
$\Delta.E_2 \geq 2$.  By Lemma \ref{A} there is an $A > 0$ primitive such that $A^2 = 0$ and $2E + E_1 - E_2
\sim A + k \Delta$, where $k:= - \Delta.(2E + E_1 - E_2) \geq 2$. From $4 = E_2.(2E + E_1 -E_2) = E_2.A + k
E_2.\Delta$ we find that $E_2 \eqv A$ and $k = 2$, so that $2E + E_1 - E_2 \eqv E_2 + 2\Delta$. Hence $E_1 \eqv
2(E_2 + \Delta - E)$.
\end{proof}

\section{A couple of useful results using vector bundles methods}
\label{cliff}
 
In the present section we will derive two useful results from the well-known vector bundle methods
introduced by Lazarsfeld and Tyurin (\cite{gl, la, ty}). The methods will be pushed a little bit forward on an
Enriques surface. To this end, recall that if $C$ is a smooth irreducible curve on a smooth irreducible surface
$S$ with $h^1(\O_S) = 0$ and $A$ is a globally generated line bundle on $C$, one can construct
(\cite{la, cp, par}) a vector bundle $\E(C,A)$ of rank $h^0(A)$, with $\det \E(C,A) = \O_S(C)$ and fitting into
an exact sequence
\begin{eqnarray}
\label{eq:1}
0 \hpil H^0(A)^{\ast} \* \O_S \hpil \E(C,A) \hpil \N_{C/S} \* A^{-1} \hpil 0.
\end{eqnarray}

We will make use of the following variant of a well-known result in \cite{dm, kn, glmPN}:  

\begin{prop} 
\label{mainbog}
Let $|L|$ be a base-point free complete linear system on an Enriques surface $S$ and assume that there
is a smooth irreducible curve $C \in |L|$ with a base-point free line bundle $A$ such that $|A|$ is
a $g^1_k$ and $L^2 \geq \max \{4k-2, 2k - 2 + 2\phi(L) \}$. Then either
\begin{itemize}
\item[(a)] there is an effective nontrivial decomposition $L \sim N + N'$ such that $|N'|$ is  
base-component free with $2(N')^2 \leq L.N' \leq  (N')^2 + k \leq 2k$ and $N'_{|C} \geq A$.
Either $N \geq N'$ or $|N|$ is base-component free and $\Bs |N| \subset C$ and $H.(N - N') \geq 0$ for any ample
$H$. Moreover if $\phi(N') = 1$ and $L.N' \geq (N')^2 + k - 1$, then $\Bs |N'| \cap C \neq \emptyset$,   

or
\item[(b)] $L^2 = 4k - 2$ and for any $\Delta \geq 0$ such that $h^0(\E(C,A)(-\Delta)) > 0$ we can find a
line bundle $N \geq \Delta$ and a nodal cycle $R$ so that $L \sim 2N + R + K_S$, $R^2 = -2$, $h^0(R) = 1$,
$h^0(R + K_S) = 0$, $|N'| : = |N + R + K_S|$ is base-component free with $\Bs |N'| \subset C$, $N.(N + R) = k$
and
$(N + R + K_S)_{|C} \geq A$.
\end{itemize}
\end{prop}

\begin{proof}
Let $E > 0$ be any divisor such that $E^2 = 0$ and $E.L = \phi(L)$. Then  $h^0(\N_{C/S} - A - E_{|C}) =
h^0(\omega_C-A-(E+K_S)_{|C}) \geq h^1(A) - \phi(L) = 2 + \frac{1}{2}L^2 - k - \phi(L) > 0$ by hypothesis. In
particular $h^0(\N_{C/S} - A) > 0$ so that $\E: = \E(C,A)$ is globally generated off a finite set. Also we have
$h^2(\E \* \omega_S) = 0$ and $c_1(\E)^2  - 4 c_2(\E) \geq - 2$, whence, as in \cite{dm, kn, glmPN}, there exist
two line bundles $N$, $N'$ on $S$ and $Z \subset S$ with $\dim Z = 0$, such that 
\begin{equation}
\label{eq:seqbog}
0 \hpil N \hpil \E \hpil \I_{Z/S} \* N' \hpil 0  
\end{equation}
with $L \sim N + N'$ and $k = N.N' + \length(Z)$.

Now case (a) corresponds to the case $c_1(\E)^2  - 4 c_2(\E) \geq 0$, using (\ref{eq:1}) and
(\ref{eq:seqbog}). Suppose that $\phi(N') = 1$ and $L.N' \geq (N')^2 + k - 1$. Then $|N'|$ has two base points
by \cite[Thm.4.4.1]{cd}. Since $N'$ is a quotient of $\E$ off $Z$ and $\E$ is globally generated outside a
finite set contained in $C$, the two base points of $|N'|$ must lie on $C \cup Z$. But $\length(Z) = (N')^2 + k
- L.N' \leq 1$, whence $\Bs |N'| \cap C \not= \emptyset$, as stated. The rest of (a) is proved similarly to 
\cite{dm}, \cite[\S 3]{kn}, \cite[Lemma2.1]{glmPN}. Case (b) corresponds to the case $c_1(\E)^2  - 4 c_2(\E) =
-2$ and is proved as in \cite[\S 3]{kn} using \cite[Thm.3.4]{kim}. 
\end{proof}

We will also need the following simple result.

\begin{lemma} 
\label{cor:nuovo}
Let $C$ be a smooth irreducible curve on an Enriques surface $S$ and let $A$ be a base-point free line 
bundle on $C$ with $h^0(A)=2$ and $\deg A=k$. Assume that there are two line bundles $M>0, N>0$ such that $C +
K_S \sim M + N$, $h^1(M + K_S) = 0$ and $\frac{1}{2}(M^2 + N^2) + 2 > k$. Then either $h^0(\E(C,A)(- M)) > 0$ 
and $(N + K_S)_{|C} \geq A$ or $h^0(\E(C,A)(- N)) > 0$ and $(M + K_S)_{|C} \geq A$.  
\end{lemma}

\begin{proof}
Set $\E = \E(C,A)$. If $h^0(\O_C(N + K_S)(- A)) = 0$ we get, using Serre duality, the dual of
\eqref{eq:1} tensored with $\O_S(M)$, Riemann-Roch and the hypotheses, that 
\begin{eqnarray*}
h^0(\E( - N)) & \geq & h^0(A) \cdot h^0(M + K_S) - h^0(\O_C(M + K_S)(A)) = \\
& = & 2\chi(M + K_S) - \chi(\O_C(M + K_S)(A)) = \frac{1}{2}(M^2 + N^2) + 2 - k > 0.
\end{eqnarray*}
Similarly, if  $h^0(\O_C(N + K_S)(- A)) > 0$, we find that $h^0(\E(- M))>0$. To conclude just tensor
(\ref{eq:1}) by $\O_S(- M)$ or $\O_S(- N)$.
\end{proof}

\section{A framework for the study of generic gonality}
\label{setting}

The goal of this section will be to devise a method to study the gonality of general
curves $C$ in a given complete linear system $|L|$ on a surface. While all previous means of
investigation  are essentially based on the vector bundle method of Bogomolov, Lazarsfeld, Tyurin and
others, and therefore work when instability conditions hold, requiring $L^2$ to be large enough, our
approach  will be, in many cases, independent of such conditions. The idea will be a sort of "liaison"
using the zero-dimensional schemes defining the gonality.

We will often use the ensuing two definitions.

\begin{defn}
\label{gengon}
Let $S$ be an Enriques surface, let $|L|$ be a base-point free complete linear system on $S$, $|L|_{sm}$
the open subset of smooth curves in $|L|$ and let $k \geq 3$ be an integer. We will say that a nonempty 
open subset $\U \subset |L|_{sm}$ has {\bf generic gonality $k$} if the following two conditions
hold:
\begin{itemize}
\item[(i)] $\gon(C) = k$ for all $C \in \U$;
\item[(ii)] if $k = 2 \phi(L)$ then $k^2 < 2L^2$ and every $C \in \U$ carries (at least) one
$g^1_k$, $A_C$, that is not cut out by a genus one pencil on $S$.
\end{itemize}
\end{defn}

\begin{defn}
\label{dominant}
Let $(L, k)$ be as in Definition {\rm \ref{gengon}} and let $\U \subset |L|_{sm}$ be a nonempty open
subset with generic gonality $k$. For every $C \in \U$ let $A_C$ be a $g^1_k$ on $C$ such that, 
if $k = 2\phi(L)$, $A_C$ is not cut out by a genus one pencil on $S$. 
Let $\{D_1, \ldots, D_n\}$ be a set of divisors on $S$. 

We will say that $\{D_1, \ldots, D_n\}$ is a {\bf $(\U, L, k, \{A_C\}_{C \in \U})$-dominant
set of divisors} if for every $C \in \U$ there is an $i \in \{1, \ldots, n \}$ 
such that $|D_i|$ is base-component free, $L - D_i > 0$, $h^1(D_i - L) = 0$ and $(D_i)_{|C} \geq A_C$.
\end{defn}

Let us outline our strategy to study the generic gonality of curves in a complete linear
system $|L|$.  To prove that a general curve $C \in |L|$ has gonality $k + 1$, assume that there is a
nonempty open  subset $\U \subset |L|_{sm}$ of generic gonality $k$. Make a choice of a $g^1_k$ for
every $C \in \U$  and suppose to have found a $(\U, L, k, \{A_C\}_{C \in \U})$-dominant set of divisors
$\{D_1, \ldots, D_n\}$. Intersecting them with $C$ we find some ``unexpected" subschemes $Z_h$ and this
will imply the existence of a line bundle $M$ that will cut out, outside of $Z_h$, the given $g^1_k$ on
every curve $C \in \U$.  This will lead to some inequalities on $L^2, k, L.M$ etc. (Proposition
\ref{lemma:somma}) that will be used either to deduce a contradiction or to prove that the gonality is
computed by $\mu(L)$.

The crux of the argument will be the following construction, with the help of the $Z_h$'s,
of a suitable incidence correspondence.

Let $(L, k)$ be as in Definition {\rm \ref{gengon}} and let $\U \subset |L|_{sm}$ be a
nonempty open subset with generic gonality $k$. Let $M \geq 0$ be an effective line bundle on $S$
with $M^2 \geq 4$. Define
\[ h_M = M.L - k, \]
and, for every zero-dimensional subscheme $X \subset S$ of length $\frac{1}{2}M^2 - 1$, define
\begin{equation}
\label{w}
W(X) = \Bs |\I_{X/S} \* M|
\end{equation}
and the following incidence subvariety of $\Hilb^{\frac{1}{2}M^2-1}(S) \times |L|$:
\begin{equation}
\label{jm}
\J_M = \{ (X,C) \; : \; C \in \U, X \subset C, \dim W(X) = 0, \deg (W(X) \cap C) = 
h_M, h^0(\I_{X/S} \* M) = 2 \}
\end{equation}
(where $W(X) \cap C$ denotes the scheme-theoretic intersection), together with its two projections $\pi_{1, M}
: \J_M \to \Hilb^{\frac{1}{2}M^2-1}(S)$ and $\pi_{2, M} : \J_M \to |L|$.

With these in mind we have
\begin{lemma}  
\label{oneM}
Let $(L, k)$ be as in Definition {\rm \ref{gengon}} and with $L^2 \geq 8$. Let $\U \subset |L|_{sm}$ be
a nonempty open subset with generic gonality $k$ and let $\{D_1, \ldots, D_n\}$ be a 
$(\U, L, k, \{A_C\}_{C \in \U})$-dominant set of divisors. Then there exists a nonempty subset
$\{D_{i_1}, \ldots, D_{i_s}\}$ of $\{D_1, \ldots, D_n\}$ and, for every $j = 1, \ldots, s$, a 
base-point free line bundle $M_{i_j}$ on $S$ such that:
\begin{itemize}
\item[(i)] $0 \leq M_{i_j} \leq D_{i_j}$;
\item[(ii)] $4 \leq M_{i_j}^2 \leq D_{i_j}^2$ and $M_{i_j} = D_{i_j}$ if $M_{i_j}^2 = D_{i_j}^2$;
\item[(iii)] $M_{i_j}.L \leq k + M_{i_j}^2$.
\item[(iv)] $\phi(D_{i_j}) \geq \phi(M_{i_j}) \geq 2$.
\end{itemize}
Furthermore, $\pi_{2, M_{i_j}}$ is a dominant morphism for at least one $j$.
\end{lemma}
\noindent {\it Proof}.
Given $C \in \U$ we know by Definition \ref{dominant} that there is an $i \in \{1, \ldots, n \}$ 
such that $|D_i|$ is base-component free, $L - D_i > 0$, $h^1(D_i - L) = 0$ and $(D_i)_{|C} \geq A_C$.
Pick $Z \in |A_C|$. Then $h^0(\I_{Z/S} \* D_i) = h^0((D_i)_{|C} - A_C) \geq 1$, 
whence there is a divisor $D_{i, Z} \in |\I_{Z/S} \* {D_i}|$ and we can write $D_{i, Z} \cap C = Z + T(C, Z,
D_{i, Z})$ as divisors on $C$. We record for later use that $\deg(T(C, Z, D_{i, Z})) = D_i.L - k$. Set, for
simplicity, $T_i = T(C, Z, D_{i, Z})$. Then 
\[ 0 \hpil D_i - L \hpil \I_{T_i/S} \* D_i \hpil A_C \hpil 0 \] 
gives that $h^0(\I_{T_i/S} \* D_i) = 2$. Defining $B_i = B_{Z, C, D_{i, Z}}$ to be the (possibly empty) 
base-divisor of $|\I_{T_i/S} \* D_i|$ and $M_i = M_{Z, C, D_{i, Z}} := D_i - B_i$, we can write
\begin{equation}
\label{subpencil}
|\I_{T_i/S} \* D_i| = |V_i| + B_i,
\end{equation}
where $V_i = V_{Z, C, D_{i, Z}} \subseteq H^0(M_i)$ is a pencil without base components. In particular 
$\dim |M_i| \geq 1$ and $|M_i|$ is base-component free. Note also that any such $M_i$ satisfies 
$0 \leq M_i \leq D_i$. Since there are only finitely many ways of decomposing $D_i$ into the sum of two
effective divisor classes, we see that, in any case, the possible line bundles $M_i = M_{Z, C, D_{i,
Z}}$ are  finitely many. 

\begin{sublemma} 
\label{lemma:nobase}
For any $C \in \U$, $Z \in |A_C|$, $D_{i, Z} \in |\I_{Z/S} \* D_i|$ fix $B_i = B_{Z, C, D_{i, Z}}$, 
$M_i = M_{Z, C, D_{i, Z}}$ and $V_i = V_{Z, C, D_{i, Z}}$ as above. Then there are two divisors $T_{M_i}
=  T_{M_i}(C, Z, D_{i, Z})$, $T_{B_i} = T_{B_i}(C, Z, D_{i, Z})$ giving an effective decomposition 
$T_i = T_{M_i} + T_{B_i}$ such that
\[ |\I_{T_i/S} \* D_i| = |\I_{T_{M_i}/S} \* M_i| + B_i, \ H^0(\I_{T_{M_i}/S} \* M_i) = V_i \] 
and with $T_{B_i} = C \cap B_i$, ${M_i}_{|C} \sim A_C + T_{M_i}$, $h_{M_i} = \deg T_{M_i} \leq M_i^2$. 
Moreover (ii) of Lemma {\rm \ref{oneM}} holds and, possibly after shrinking $\U$ to some nonempty open
subset $\U'$, we can assume that no $C \in \U'$ passes through any of the two base points of all the
possible $|M_j|$ that have $\phi(M_j) = 1$.
\end{sublemma}
\begin{proof}
For simplicity set $B = B_i, M = M_i, D = D_i, T = T_i, V = V_i$. If $B = 0$ we just choose $T_M = T$ and $T_B
= 0$. Suppose now that $B > 0$. Set $T_M := \Bs |V_{|C}|$. Then  $T_M \leq T$ as $A_C$ is base-point
free. Set $T_B := T - T_M$. Then one easily sees that $T_B \leq C \cap B$. By definition of $T_M$ we have
an inclusion $V \subseteq H^0(\I_{T_M/S} \* M)$ and $h^0(\I_{T_M/S} \* M) \geq 2$. Therefore
$h^0(M_{|C} - T_M) \geq 2$. But $\deg (M_{|C} - T_M)  \leq k$, whence $T_B = C \cap B$, as stated, using the
assumption $\gon(C) = k$. Now $B_{|C} + M_{|C} \sim A_C + T_M + T_B$ and it follows that $M_{|C} \sim A_C +
T_M$ and $\deg T_M = M.L - k = h_M$. This implies that $h^0(M_{|C} - T_M)= h^0(\I_{T_M/S} \* M) = 2$, whence $V
= H^0(\I_{T_M/S} \* M)$ and, by (\ref{subpencil}), $|\I_{T/S} \* D| = |\I_{T_M/S} \* M| + B$.

Since there are only finitely many $M$'s, we can shrink $\U$ to a nonempty open subset 
$\U'$ such that no $C \in \U'$ passes through any of the two base points of all the possible $|M_j|$ 
that have $\phi(M_j) = 1$. 

Now two general distinct elements $M'$ and $M''$ in $|V| = |\I_{T_M/S} \* M|$ have no
common components, so that $W: = M' \cap M''$ has dimension zero and degree $M^2$. As $T_M \subseteq
W$ we deduce that $\deg T_M \leq M^2$. If $M^2 = 0$ we get that $M_{|C} \sim A_C$, contradicting
(ii) of Definition \ref{gengon}.

If $M^2 = 2$ then $|\I_{T_M/S} \* M| = |M|$ and $|M|$ has two distinct base points $x$
and $y$. Now $T_M \subseteq W \cap C = \Bs |M| \cap C = \emptyset$.
Therefore $T_M = 0$ and $M_{|C} \sim A_C$, whence $k = M.C \geq 2\phi(L)$ (since $h^0(M)
= 2$). As $\gon(C) = k$ we also have $k \leq 2\phi(L)$, therefore we must have $k = 2\phi(L)$,
whence, by (ii) of Definition \ref{gengon}, $(M.L)^2 = k^2 < 2L^2 = M^2 L^2$, which
contradicts the Hodge index theorem.  

Therefore $M^2 \geq 4$. Now, as $D$ and $M$ are nef, we get $D^2 = D.(M + B) \geq D.M = (M + B).M
\geq M^2$. If equality holds, then $H^0(M) = H^0(D)$ by Riemann-Roch, so that $B$ is a base-component of $|D|$,
whence $B = 0$  and (ii) of Lemma \ref{oneM} is proved.
\end{proof}

\noindent {\it Continuation of the proof of Lemma {\rm \ref{oneM}}}.
Let $M_i$ and $D_i$ be the line bundles arising in Sublemma \ref{lemma:nobase}.
If, on any $C \in \U'$, we choose a divisor $Z \in |A_C|$ and a divisor $D_{i, Z} \in |\I_{Z/S} \* D_i|$, it is
easily seen that we can find a divisor $X_i \leq T_{M_i}$ such that
\[ 2 = h^0(\I_{T_{M_i}/S} \* M_i) = h^0(\I_{X_i/S} \* M_i) = h^0(M_i) - \deg X_i = 
\frac{1}{2} M_i^2 + 1 - \deg X_i, \]
therefore with $h^0(\I_{X_i/S} \* M_i) = 2$ and $\deg X_i = \frac{1}{2} M_i^2 - 1$. Hence 
$X_i \in \Hilb^{\frac{1}{2} M_i^2 - 1}(S)$ and, as in (\ref{w}), $W(X_i) = \Bs |\I_{X_i/S} \* M_i|$ 
is such that $\dim W(X_i) = 0$, $\length(W(X_i)) = M_i^2$ and $W(X_i) = M_i' \cap M_i''$ for two 
general divisors $M_i', M_i'' \in |\I_{X_i/S} \* M_i|$. 

Note that also $T_{M_i} \subseteq W(X_i) \cap C$, the scheme-theoretic intersection. In fact we have
\begin{claim} 
\label{cl:tutto}
$T_{M_i} = W(X_i) \cap C$ and $\phi(D_i) \geq \phi(M_i) \geq 2$. 
\end{claim}

\begin{proof}
Set $Z' = W(X_i) \cap C$ and suppose that $\deg Z' > \deg T_{M_i}$. Then 
$h^0({M_i}_{|C} - Z') \geq h^0(\I_{Z'/S} \* M_i) \geq h^0(\I_{W(X_i)/S} \* M_i) = 2$. 
On the other hand $\deg ({M_i}_{|C} - Z') < M_i.L-h_{M_i}=k$, contradicting our assumption that
$\gon(C) = k$.

To see the second assertion note that, since $D_i$ is nef, there is a nef divisor $E$ with $E.D_i =
\phi(D_i)$, whence $\phi(D_i) \geq \phi(M_i)$. Now assume that $\phi(M_i) = 1$. Then $W(X_i)$ contains the two
base points  of $|M_i|$. Since $W(X_i) \cap C = T_{M_i}$ and, by Sublemma \ref{lemma:nobase}, $C$ does not
contain any of them, we must have $h_{M_i} = M_i.L - k = \deg T_{M_i} \leq \length(W(X_i)) - 2 \leq M_i^2 -
2$. By Lemma \ref{lemma:nefred} we can write $M_i \sim \frac{1}{2} M_i^2 E_1 + E_2$ for $E_p > 0$ with
$E_p^2 = 0$, $p = 1, 2$ and $E_1.E_2 = 1$. Hence $(\frac{1}{2} M_i^2 + 1) \phi(L) \leq M_i.L \leq k + M_i^2 - 2
\leq 2 \phi(L) + M_i^2 - 2$. Since $\phi(L) \geq 2$, we deduce that $\phi(L) = 2, k = 2 \phi(L) = 4$ and
$E_1.L = E_2.L = 2$. Let $b := \lfloor \frac{L^2}{4} \rfloor$. Then $(L - b E_1)^2 = L^2 - 4b \geq 0$, therefore
$L - b E_1 > 0$ by Lemma \ref{lemma:nefred} and $2 - b = E_2.(L - b E_1) \geq 0$, by Lemma \ref{lemma10}.
Since $L^2 \geq 8$ by hypothesis we deduce that $b = 2$ and $E_2.(L - 2 E_1) = 0$ whence, again by
Lemma \ref{lemma10}, we have that $(L - 2 E_1)^2 = 0$, giving $L^2 = 8$. Therefore $k = 2
\phi(L)$ and $k^2 = 2L^2$, contradicting (ii) of Definition \ref{gengon}. 
\end{proof}

\renewcommand{\proofname}{Conclusion of the proof of Lemma {\rm \ref{oneM}}}  
\begin{proof}
We have therefore found a nonempty open subset $\U' \subset |L|_{sm}$ such that, for any $C \in \U'$, 
by choosing $Z \in |A_C|$ and $D_{i, Z} \in |\I_{Z/S} \* D_i|$, we have a line bundle $M_i$ satisfying 
(i), (ii), (iii) and (iv) of Lemma \ref{oneM} (and base-point free by \cite[Prop.3.1.4 and
Thm.4.4.1]{cd}) and an element $X_i \in \Hilb^{\frac{1}{2} M_i^2 - 1}(S)$  such that
\[ X_i \subset C, \ \dim W(X_i) = 0, \ \deg (W(X_i) \cap C) = h_{M_i} \ \mbox{and} \ 
h^0(\I_{X_i/S} \* M_i) = 2 \]
whence $(X_i, C) \in \J_{M_i}$, where $\J_{M_i}$ is as in (\ref{jm}). Now if
$\{i_1, \ldots, i_s \} = \{i : 1 \leq i \leq n \ \mbox{and there exists} \ C \in \U' \ 
\mbox{such that} \ (D_i)_{|C} \geq A_C \}$, we have proved that $\U' = \bigcup\limits_{j = 1}^s \Im \pi_{2,
M_{i_j}}$. Hence there must be at least one $M_{i_j}$ such that $\pi_{2, M_{i_j}}$ is a dominant morphism.
\end{proof}
\renewcommand{\proofname}{Proof}

\begin{prop}  
\label{lemma:somma}
Let $S$ be an Enriques surface and let $|L|$ be a base-point free complete linear system on $S$ with
$L^2 \geq 8$ and let $k \geq 3$ be an integer. Let $\U \subset |L|_{sm}$ be a nonempty open subset with
generic gonality $k$ and let $\{D_1, \ldots, D_n\}$ be a $(\U, L, k, \{A_C\}_{C \in \U})$-dominant 
set of divisors. Then there exists a base-point free line bundle $M$ on $S$ such that: 

\vskip .1cm

\noindent {\rm (a)} For at least one $i$ we have $M \leq D_i$, $|D_i|$ is base-component free, $L - D_i > 0$, 
$h^1(D_i - L) = 0$, $\phi(D_i) \geq \phi(M) \geq 2$, $4 \leq M^2 \leq D_i^2$ and $M = D_i$ if 
$M^2 = D_i^2$. Moreover $M.L \leq k + M^2$.

\noindent {\rm  (b)} If $h_M := M.L - k \geq M^2 - 1$ we have
\begin{equation}
\label{eq:databe}
3M.L \leq 3M^2 + \frac{1}{2} L^2 + 2h^1(L - M) + h^1(L - 2M) + k - 1
\end{equation}
and
\begin{equation}
\label{eq:datad}
M^2 - M.L + k + 2h^1(L - M) + h^0(2M - L + K_S) \geq 2. 
\end{equation}

\noindent {\rm (c)} If equality occurs in {\rm (\ref{eq:databe})} then
\begin{equation}
\label{eq:datac}
M.L = M^2 + k \ \mbox{and} \ h^0(L - 2M) = 0.   
\end{equation}
\end{prop} 

\begin{proof}
Let $M = M_i$ be a line bundle arising from Lemma \ref{oneM} for which $\pi_{2, M_i}$ is a dominant
morphism and set $D = D_i$. Then (a) follows from Lemma \ref{oneM}.
To see (b) we study the fibers of the map $\pi_{1, M}$ of $\J_M$ in (\ref{jm}). For $X \in \Im \pi_{1, M}$
consider the set 
\[ \Sigma_X = \{ Y \subset S \; : \; Y = W(X) \cap C \;  \mbox{for some} \; C \in \U \}, \]
where $W(X) \cap C$ denotes the scheme-theoretic intersection. By definition of $\J_M$ we have, for any $Y \in
\Sigma_X$, that $\deg Y = \deg (W(X) \cap C) = h_M \geq M^2 - 1$, while $Y \subseteq W(X)$ and $\length(W(X)) =
M^2$. Therefore either $h_M = M^2$ and $Y = W(X)$ or $h_M = M^2 - 1$ and $Y$ is a maximal subscheme of $W(X)$.
Now $W(X) = \Bs |\I_{X/S} \* M|$ and $h^0(\I_{X/S} \* M) = 2$, whence $W(X)$ is a locally complete intersection
subscheme of $S$, therefore it is Gorenstein, whence it has finitely many maximal subschemes.
Therefore, in either case, we deduce that $\Sigma_X = \{Y_1, \ldots, Y_t \}$. Now if $(X, C) \in
(\pi_{1, M})^{-1}(X)$ then $W(X) \cap C = Y_j$ for some $j = 1, \ldots, t$. Therefore we have the
following inclusion in $|L|$
\[(\pi_{1, M})^{-1}(X) \sub \bigcup\limits_{j = 1}^t |H^0(\I_{Y_j/S} \* L)| \]
whence
\begin{equation}
\label{eq:a8}
\dim (\pi_{1, M})^{-1}(X) \leq \max \{h^0(\I_{Y/S} \* L) - 1, Y \in \Sigma_X \}.
\end{equation}
To estimate $h^0(\I_{Y/S} \* L)$ we will use the exact sequence
\[ 0 \hpil \I_{W(X)/S} \* L \hpil \I_{Y/S} \* L \hpil \I_{Y/W(X)} \* L \hpil 0, \] 
and 
\begin{equation}
\label{eq:a8'}
h^0(\I_{Y/W(X)} \* L) = \chi (\I_{Y/W(X)} \* L) = \chi (\O_{W(X)} \* L) - 
\chi (\O_{Y} \* L) = M^2 - h_M. 
\end{equation}
Since $W(X) = M' \cap M''$ for two distinct $M', M'' \in |\I_{X/S} \* M|$, we have
\[ 0 \hpil L - M \hpil \I_{W(X)/S} \* L \hpil \I_{W(X)/M'} \* L \hpil 0, \]
whence
\begin{eqnarray*}
h^0(\I_{W(X)/S} \* L) & \leq & h^0(L - M) + h^0(\I_{W(X)/M'} \* L) = h^0(L - M) + h^0(\O_{M'}(L - M)) \leq \\
\nonumber & \leq & 2h^0(L - M) + h^1(L - 2M).
\end{eqnarray*}
Combining the latter inequality with (\ref{eq:a8}) and (\ref{eq:a8'}) we find
\[ \dim (\pi_{1, M})^{-1}(X) \leq M^2 - h_M + 2h^0(L - M) + h^1(L - 2M) - 1. \]
Note now that the latter bound does not depend any more on $X$. Hence, using the fact that $\dim
\Hilb^{\frac{1}{2}M^2-1}(S) = M^2 - 2$, we have, for any irreducible component $\J$ of $\J_M$, that
\begin{equation}
\label{eq:a10}
\dim \J \leq 2M^2 - 3 + 2h^0(L - M) + h^1(L - 2M) - h_M.
\end{equation}
Since $\pi_{2,M}$ is a dominant morphism, there must be a component $\J_0$ of $\J_M$ such that 
$\dim \J_0 \geq \dim |L|= \frac{1}{2}L^2$, whence
\begin{equation}
\label{vecchia34}
\frac{1}{2}L^2 \leq 2M^2 - 3 + 2h^0(L - M) + h^1(L - 2M) - M.L + k.
\end{equation}
Since $h^2(L - M) \leq h^2(L - D) = 0$, we obtain (\ref{eq:databe}) combining Riemann-Roch and
(\ref{vecchia34}).

On the other hand the inclusions $Y \subseteq W(X) = M' \cap M'' \subset M'$ give the
exact sequences
\[ 0 \hpil L - M \hpil \I_{Y/S} \* L \hpil \I_{Y/M'} \* L \hpil 0, \]
\[ 0 \hpil \O_{M'}(L - M) \hpil \I_{Y/M'} \* L \hpil \I_{Y/W(X)} \* L \hpil 0. \]
Therefore 
\begin{equation}
\label{eq:a12'''}
h^1(\I_{Y/M'} \* L) \leq h^1(\O_{M'}(L - M)) \leq h^1(L - M) + h^2(L - 2M)   
\end{equation}
and 
\[ h^1(\I_{Y/S} \* L) \leq h^1(L - M) + h^1(\I_{Y/M'} \* L) \leq 2 h^1(L - M) + h^2(L - 2M), \]
whence
\[ h^0(\I_{Y/S} \* L) \leq h^0(L) - h_M + 2h^1(L - M) + h^2(L - 2M), \]
so that, as above,
\begin{eqnarray*}
\dim \J \leq  \dim |L| + M^2 - 2 - h_M + 2h^1(L - M) + h^0(2M - L + K_S)
\end{eqnarray*}
for every irreducible component $\J$ of $\J_M$. Since $\pi_{2, M}$ is a dominant morphism, there must 
be a component $\J_0$ of $\J_M$ such that $\dim \J_0 \geq \dim |L|$, whence (\ref{eq:datad}) 
follows. 

Assume now that equality occurs in (\ref{eq:databe}). Then, by Riemann-Roch, 
$h^0(L - 2M) + h^0(2M - L + K_S)  = h_M - M^2 + 2 - 2h^1(L - M)$.
Recalling that $M^2 - 1 \leq h_M \leq M^2$, we get the following three possibilities, where we set,
for simplicity, $\alpha = h^0(L - 2M) + h^0(2M - L + K_S)$ and $\beta = h^1(L - M)$:
\begin{equation}
\label{eq:four}
(\alpha, h_M, \beta) = (0, M^2, 1), (2, M^2, 0), (1, M^2 - 1, 0). 
\end{equation}
Now equality must also occur in (\ref{eq:a10}) for $\J = \J_0$, so that it follows that
the restriction of $\pi_{1, M}$ to $\J_0$ is a dominant morphism. We know that there is a nonempty 
open subset $\V \subset \Hilb^{\frac{1}{2}M^2-1}(S)$ of zero-dimensional subschemes that impose
independent conditions to $|M|$, whence the incidence correspondence 
\[ \I_M = \{ (X, C) \; : \; X \in \V, C \in |M|, X \subset C \} \]
has both projections to $\V$ and to $|M|$ that are dominant morphisms. Therefore a general element 
$X \in \Hilb^{\frac{1}{2}M^2-1}(S)$ lies on a general element $M' \in |M|$. In particular we can pick 
$X \in \Hilb^{\frac{1}{2}M^2-1}(S)$ general so that $W(X) = M' \cap M''$ with $M'$ a smooth
irreducible curve in $|\I_{X/S} \* M|$. Moreover, in case $h^0(2M - L + K_S) = 1$, if we let
$\Delta$ be the unique element of $|2M - L + K_S|$, we can also assume that a general element $X \in
\Hilb^{\frac{1}{2}M^2-1}(S)$ does not intersect $\Delta$. Now let $\varphi_M : S \to
\PP^{\frac{1}{2}M^2}$ be the morphism defined by $|M|$. Since $X \cap \Delta = \emptyset$ and
$\varphi_M(S)$ is a surface (since $M^2 \geq 4$), we can further assume that the curve (or point)
$\varphi_M(\Delta)$ does not intersect a general $(\PP^{\frac{1}{2}M^2-2})$-plane in
$\PP^{\frac{1}{2}M^2}$.

Now, in the case $(\alpha, h_M, \beta) = (1, M^2 - 1, 0)$, there is a point $p \in S$ such that 
\[ Y + p = W(X) \; \mbox{as divisors on} \; M'. \]
We henceforth set $W' = W(X) - Y$ on $M'$, so that $W' = p$ in the case $(\alpha, h_M, \beta) =
(1, M^2 - 1, 0)$ and $W' = 0$ in the two other cases in (\ref{eq:four}). Moreover, since
$M'$ is smooth we have $\I_{Y/M'} \cong \I_{W(X)/M'}(W') \cong \O_{M'}(-M)(W')$. 

We have $M^2 - h_M + 2\beta + \alpha - 2 = 0$ in all the cases in (\ref{eq:four}), and
since $h^0(2M - L + K_S) \leq \alpha$, we must have by (\ref{eq:datad}) that $\alpha = h^0(2M -
L + K_S)$, $h^0(L - 2M) = 0$ and therefore that equality occurs in (\ref{eq:datad}), whence it must also
occur in (\ref{eq:a12'''}), so that 
\begin{equation}
\label{eq:a13'}
h^1(\O_{M'}(L - M)(W')) = h^1(L - M) + h^2(L - 2M) = \beta + \alpha. 
\end{equation}
In case $(\alpha, h_M, \beta) = (1, M^2 - 1, 0)$ we have $W' = p$ and $h^1(M - L + K_S) = \beta =
0$ and from (\ref{eq:a13'}), Serre duality and
\[ 0 \hpil M - L + K_S \hpil \I_{\{p\}/S} \* (2M - L + K_S) \hpil \O_{M'}(2M - L + K_S)(- p) \hpil
0 \] 
we see that
\[ h^0(\I_{\{p\}/S} \* (2M - L + K_S)) \geq h^0(\O_{M'}(2M - L + K_S)(- p)) = h^1(\O_{M'}(L -
M)(p)) = 1 \]
whence $h^0(2M - L + K_S) = h^0(\I_{\{p\}/S} \* (2M - L + K_S)) = 1$ . This means that the
unique element $\Delta$ of $|2M - L + K_S|$ contains $p$. Since $X \cap \Delta = \emptyset$ we have
$p \not \in X$, but since $p \in W(X)$ and $H^0(\I_{W(X)/S} \* M) = H^0(\I_{X/S} \* M)$ we have that
$\varphi_M(p)$ is contained in the linear span of $\varphi_M(X)$, which is a general
$\PP^{\frac{1}{2}M^2-2}$. Since $\varphi_M(p) \in \varphi_M(\Delta)$, this is a contradiction.

Therefore the third case in (\ref{eq:four}) does not occur. Hence $h^0(L - 2M) = 0$,
$M.L = M^2 - k$ and (\ref{eq:datac}) is proved.
\end{proof}

Now we need a simple, auxiliary result:

\begin{lemma} 
\label{lemma:A}
Let $L$ be a line bundle on an Enriques surface $S$ with $L^2 > 0$ and assume $L \sim M + N$ 
is a nontrivial, effective decomposition such that $M^2 \geq 0$. Set $l = M.N$.
\begin{itemize}
\item[(a)] If $l \leq 2\phi(L) - 2$, $M^2 \geq 4$, $h^0(2M - L + K_S) \leq 2$ and $L^2 \geq 4l + 4$, 
then there is a line bundle $M_0 > 0$ such that $M_0^2 = 4$ and $M_0.L \leq l + 4$.
\item[(b)] If $L$ is nef, $N^2 \geq 0$, $l \leq 2\phi(L) - 1$, $L^2 \geq 4l - 4$ and $h^1(M) \neq 0$, 
then there is a line bundle $M_0 > 0$ such that $M_0^2 = 4$ and $M_0.L \leq l + 2$.
\end{itemize}
\end{lemma}

\begin{proof}
(a) We have $(L - 2M)^2 = (M - N)^2 = L^2 - 4M.N = L^2 - 4l \geq 4$ by assumption, whence by 
Riemann-Roch and our assumptions, we get that $h^0(L - 2M) \geq 3$.

If $M^2 = 4$, then $M.L = M^2 + M.N = l + 4$ and we are done.

If $M^2 \geq 6$, then pick $E > 0$ such that $E^2 = 0$ and $E.M = \phi(M)$. Then $(M - E)^2 
\geq 2$ whence $M - E > 0$ by Lemma \ref{lemma:nefred}. Moreover, $(M - E).(N + E) \leq l - 1$, 
since $E.(L - 2M) \geq \phi(L - 2M) \geq 1$. If $(M - E)^2 = 2$ we get
$2\phi(L) \leq (M - E).L \leq l + 1 \leq 2\phi(L) - 1$, a contradiction. Hence $(M - E)^2 \geq 4$, and since
$(M - E)^2 < M^2$, $h^0(2(M - E) - L + K_S) \leq h^0(2M - L + K_S) \leq 2$, $l - 1 \leq 2\phi(L) - 3$ and $L^2
> 4(l - 1)  + 4$, we can repeat the process if necessary, which must eventually end.

(b) If $M^2 = 0$, then $h^0(M) \geq 2$ implies $2 \phi(L) \leq M.L = M.N = l \leq 
2\phi(L) - 1$, a contradiction. Hence $M^2 > 0$, and by Theorem \ref{cor:qnef} there is a 
$\Delta > 0$ such that $\Delta^2 = -2$ and $\Delta.M \leq -2$. Set $M' = M - \Delta$ and $N' = N +
\Delta$. Note that $M' > 0$ by Lemma \ref{A}. We have $(M')^2 \geq 4$ and $M'.N'  \leq l-2 \leq 2\phi(L) - 3$.
Since $L^2 \geq 4l - 4 = 4(l - 2) + 4$, we are done by (a) if $h^0(2M' - L + K_S) \leq 2$.

Assume therefore that $h^0(2M' - L + K_S) \geq 3$. This implies $h^0(2N' - L + K_S) = 0$. Now
$(N')^2 = (N + \Delta)^2 \geq N^2 + 2 \geq 2$ (recall that $\Delta.N \geq -\Delta.M \geq 2$ since $L$ 
is nef). If $(N')^2 = 2$ then $2\phi(L) \leq N'.L = (N')^2 + N'.M' \leq l \leq 2\phi(L) - 1$,
a contradiction.

Therefore $(N')^2 \geq 4$ and we are done again by (a). 
\end{proof}

This allows us to prove the main result of this section:

\begin{prop} 
\label{prop:C}
Let $S$ be an Enriques surface and let $|L|$ be a base-point free complete linear system on $S$ with
$L^2 \geq 8$ and let $k$ be an integer such that $3 \leq k \leq 2\phi(L) - 1$ and $L^2 \geq 4k - 4$.
Let $\U \subset |L|_{sm}$ be a nonempty open subset with generic gonality $k$ and let 
$\{D_1, \ldots, D_n\}$ be a $(\U, L, k, \{A_C\}_{C \in \U})$-dominant set of divisors such that 
$(L - D_j)^2 \geq 0$, $L^2 \geq D_j^2 + 2k - 2$ and $h^0(2D_j - L + K_S) \leq 1$ for all 
$j = 1, \ldots, n$. Then there exists a line bundle $M_0 > 0$ such that $M_0^2 = 4$ and $M_0.L \leq 
k + 2$.
\end{prop}

\begin{proof}
Let $M$ be a line bundle obtained in Proposition \ref{lemma:somma}. Then $L - M \geq L - D_i > 0$ 
for some $i$, $M.L - M^2 \leq k$ and $4 \leq M^2 \leq D_i^2$. Hence $h^0(2M - L + K_S) \leq 
h^0(2D_i - L + K_S) \leq 1$. Since $k - 2 \leq 2\phi(L) - 3$, we are done by Lemma \ref{lemma:A}(a) if 
$M.L - M^2 \leq k - 2$.

Therefore we can assume $M.L - M^2 \geq k - 1$ and we can apply Proposition 
\ref{lemma:somma}(b).

We now claim that $(L - M)^2 \geq 0$. Indeed if $M^2 = D_i^2$ then $M = D_i$ by 
Proposition \ref{lemma:somma}(a), so that $(L - M)^2 = (L - D_i)^2 \geq 0$. On the other
hand if $M^2 \leq D_i^2 - 2$ then $(L - M)^2 \geq L^2 - 2k - D_i^2 + 2 \geq 0$. Thus if $h^1(L -
M) \neq 0$ we are done again by Lemma \ref{lemma:A}(b). Therefore we can also assume $h^1(L - M) = 0$. 
From Proposition \ref{lemma:somma}(b) we get
\[ 2 \leq M^2 - M.L + k + h^0(2M - L + K_S) \leq h^0(2M - L + K_S) + 1, \]
whence $h^0(2M - L + K_S) = 1$. From Riemann-Roch and \eqref{eq:databe} we get
\[ 3M.L \leq 3M^2 + \frac{1}{2}L^2 + h^1(L - 2M) + k - 1 = 2M.L + M^2 + k - 1, \]
whence $M.L - M^2 = k - 1$ from  our assumptions, so that 
equality occurs in \eqref{eq:databe}. But this contradicts Proposition \ref{lemma:somma}(c). 
\end{proof}

Let $L$ be a nef line bundle on an Enriques surface $S$ with $L^2 \geq 8$ and $\phi(L) \geq
2$, so that \cite[Prop.3.1.6, Prop.3.1.4, Prop.4.5.1 and Thm.4.5.4]{cd} $|L|$ is base-point free and
a general curve $C_{\eta} \in |L|$ is not hyperelliptic. We henceforth set
\[\Sigma_s = \{C \in |L|_{sm} : \gon(C) \leq s \}. \]
As is well-known $\Sigma_s$ is a closed subset of $|L|_{sm}$.

The following two results, direct applications of Proposition \ref{prop:C}, will also be 
key results to prove Theorem \ref{main}. In their proofs we will use the following simple

\begin{rem} 
\label{usefulrem}
Let $L$ be a line bundle on an Enriques surface $S$ such that $L^2 > 0$ and let $b$ be an integer such
that $L^2 \geq 4b + 2$. Then $L^2 \geq 2b + 2\phi(L)$.
\end{rem}

\begin{lemma} 
\label{lemma:D}
Let $S$ be an Enriques surface and let $|L|$ be a base-point free complete linear system on $S$. Assume
that  a general smooth curve in $|L|$ has gonality $k \geq 3$ and that either
\begin{itemize}
\item[(i)] $L^2 \geq 4k$, or
\item[(ii)] $L^2 = 4k - 2$ and there is an effective nontrivial decomposition $L \sim M + N$
such that $M^2 \geq 0$, $N^2 \geq 0$ and $k - 1 \leq M.N \leq k$.
\end{itemize}
Then $k \geq \min \{2\phi(L), \mu(L) \}$. 
\end{lemma}

\begin{proof}
As $|L|$ is base-point free we have $\phi(L) \geq 2$ by \cite[Prop.3.1.6, Prop.3.1.4 and Thm.4.4.1]{cd}. We 
assume $k \leq 2\phi(L) - 1$ and we show the existence of a line bundle $M_0 > 0$ such that $M_0^2 = 4$ and
$M_0.L \leq k + 2$. This will be enough since if $\phi(M_0) = 1$ then by Lemma \ref{lemma:nefred} we can write
$M_0 \sim 2F_1 + F_2$ with $F_i > 0, F_i^2 = 0$ whence the contradiction $3 \phi(L) \leq M_0.L \leq k + 2 \leq
2\phi(L) + 1$.

By Remark \ref{usefulrem} with $b = k - 1$ we see that we can apply Proposition \ref{mainbog}.

Let $C \in |L|$ be a smooth curve of gonality $k$ and let $|A_C|$ be a $g^1_k$ on $C$.
Let $L \sim N_C + N'_C$ be the decomposition as in Proposition \ref{mainbog}. Then 
$\O_C(N'_C) \geq A_C$ and $|N'_C|$ is base-component free, in particular $(N'_C)^2 \geq 0$.

We now show that we can assume that 
\begin{equation} \label{eq:inter}
N_C^2 \geq 0, \; L^2 \geq 2k + (N'_C)^2 \; \mbox{and} \; h^0(2N'_C - L + K_S) \leq 1.
\end{equation} 
Indeed, if we are in case (a) of Proposition \ref{mainbog}, then either $N_C.L \geq N'_C.L$, whence $N_C^2 \geq
(N'_C)^2 \geq 0$, or $|N_C|$ is base component-free, whence again $N_C^2 \geq 0$. We also have that
$(N'_C)^2 \leq k$ whence $L^2 \geq 4k - 2 \geq 2k + (N'_C)^2$. Finally, $H.(N_C - N'_C) \geq 0$ for an 
ample $H$, whence $h^0(2N'_C - L + K_S) \leq 1$.

If we are in case (b) of Proposition \ref{mainbog}, then $L^2 = 4k - 2$, so that by our assumptions there is a
decomposition $L \sim M + N$ as in (ii). By Lemma \ref{lemma:A}(b) we can assume that $h^1(M + K_S) = h^1(N) =
0$. Since
\[ \frac{1}{2}(M^2 + (N + K_S)^2) + 2 = \frac{1}{2}(L^2 - 2M.N) + 2 \geq k + 1,\]
we have that either $h^0(\E(C,A_C)(- M)) > 0$ or $h^0(\E(C,A_C)(- N + K_S)) > 0$ by Lemma \ref{cor:nuovo}.
It follows by Proposition \ref{mainbog}(b) that we can assume that either 
$N_C \geq M$ or $N_C \geq N + K_S$. In the first case we get
\[ k - 1 \leq M.N \leq M.L \leq N_C.L = N_C^2 + N_C.N'_C = N_C^2 + k \]
showing that $N_C^2 \geq 0$. Similarly we get $N_C^2 \geq 0$ also when $N_C \geq N + K_S$. Moreover, we
have that $N_C.N'_C = k$ whence $L^2 = N_C^2 + (N'_C)^2 + 2k \geq 2k + (N'_C)^2$. Finally, we note
that $h^0(2N'_C - L + K_S) = h^0(R) = 1$. 

We have therefore proved that we can assume \eqref{eq:inter}. In particular, by Lemma \ref{lemma:A}(b) we see
that we are done if $h^1(N_C + K_S) \neq 0$. Therefore we can assume that $h^1(N_C + K_S) = 0$.

Since there are only finitely many effective decompositions of $L$ we see that if 
$\U = |L|_{sm} - \Sigma_{k - 1}$, we can find a $(\U, L, k, \{A_C\}_{C \in \U})$-dominant set of 
divisors $\{D_1, \ldots, D_n\}$ all of the form $D_j = N'_C$ for some $C \in \U$. Also we have 
$h^0(2D_j - L + K_S) \leq 1$, $(L - D_j)^2 \geq 0$ and $L^2 \geq 2k + D_j^2$ for all $j$ by
\eqref{eq:inter}, whence we are done by Proposition \ref{prop:C}.
\end{proof}

\begin{lemma} 
\label{lemma:E}
Let $S$ be an Enriques surface and let $|L|$ be a base-point free complete linear system on $S$. Suppose
that  a general smooth curve in $|L|$ has gonality $k \geq 3$ with $L^2 = 4k - 4$. Assume that there 
is an effective nontrivial decomposition $L \sim M + N$ such that $0 \leq M^2 \leq 2k - 2$, 
$0 \leq N^2 \leq 2k - 2$, $M.N \leq k - 1$, $h^0(2M - L + K_S) \leq 1$ and  $h^0(2N - L + K_S) \leq 1$.
Then $k \geq \min \{2\phi(L), \mu(L) \}$. 
\end{lemma}

\begin{proof}
As in the proof of the previous lemma we can assume $k \leq 2\phi(L) - 1$ and
$h^1(M) = h^1(M + K_S) = h^1(N) = h^1(N + K_S) = 0$. 

We will prove the existence of a line bundle $B > 0$ such that $B^2 = 4$ and $B.L \leq k + 2$.

Let $\U = |L|_{sm} - \Sigma_{k - 1}$ and, for every $C \in \U$, let $|A_C|$ be a $g^1_k$ on $C$. Note that $\U$
has generic gonality $k$. As in the previous proof, by Lemma \ref{cor:nuovo}, we get that either
$N_{|C} \geq A_C$ or $(M+K_S)_{|C} \geq A_C$. 

\begin{claim}
\label{claim1}
Suppose that $N_{|C} \geq A_C$. Then $h^0(N) \geq 2$. Let $N_0$ be the moving part of $|N|$.
Then $L - N_0 > 0$, $H^1(N_0 - L) = 0$, $N_0^2 =  N^2 > 0$ and ${N_0}_{|C} \geq A_C$. 
\end{claim}

\begin{proof} We have $h^0(N) = h^0(N_{|C}) \geq h^0(A_C) = 2$. Write $N \sim \Delta + N_0$ for some $\Delta
\geq 0$. Obviously $L-N_0>0$. One relatively easily checks that ${N_0}_{|C} \geq A_C$, and using $h^1(N-L)=0$,
that the map $H^0(N_0) \to H^0({N_0}_{|C})$ is surjective, whence an isomorphism. Now if $N_0^2 > 0$ then
$h^1(N_0) = 0$, whence $H^1(N_0 - L) = 0$ and $N_0^2 = N^2 > 0$. If $N_0^2 = 0$, then $h^0({N_0}_{|C}) =
h^0(N_0) = 2$ by \cite[proof of Cor.3.1.2]{cd}, therefore ${N_0}_{|C} \sim A_C$. But then we get the
contradiction $2 \phi(L) \leq N_0.L = k$.
\end{proof}

\noindent {\it Conclusion of the proof of Lemma {\rm \ref{lemma:E}}}. 
Let $M_0$ be the moving part of $|M + K_S|$. Since a claim similar to Claim \ref{claim1} holds if
$(M + K_S)_{|C} \geq A_C$, we deduce that the pair $\{ M_0, N_0\}$ is a $(\U, L, k, \{A_C\}_{C \in
\U})$-dominant set of divisors. Moreover $h^0(2M_0 - L + K_S) \leq h^0(2M - L + K_S) \leq 1$ and $h^0(2N_0 - L
+ K_S) \leq  h^0(2N - L + K_S) \leq 1$ by assumption. To finish the proof, by Proposition \ref{prop:C}, we need
to show that $(L - N_0)^2 \geq 0$ and $N_0^2 + 2k - 2 \leq L^2$ (and similarly for $M_0$). 

By Claim \ref{claim1} we get $N_0^2 + 2k - 2 = N^2 + 2k - 2 \leq 4k - 4 = L^2$ and 
$N_0.(L - N_0) = N_0.L - N^2 \leq N.L - N^2 = N.M \leq k - 1$, from which $(L - N_0)^2 \geq 0$ easily follows. 
\end{proof}

We conclude this section by giving a few applications of the results obtained. 
They will be used to prove some cases of Theorem \ref{main}, to complete a theorem about plane
curves in \cite{kl2} and they will also be needed for the study of Gaussian maps in \cite{klGM}. 

\subsection{The strategy in low genus} 
\label{lowstr}

\hskip 1cm

We investigate here the generic gonality of some line bundles of small genus. The strategy 
that we will employ is as follows. 

Recall that $\Sigma_s = \{C \in |L|_{sm} : \gon(C) \leq s \}$. To prove that for a general 
curve $C_{\eta} \in |L|$ we have that $\gon(C_{\eta}) \geq k + 1$, we will assume that $\Sigma_k =
|L|_{sm}$ and derive a contradiction. To this end we will consider the nonempty open subset $\U =
|L|_{sm} - \Sigma_{k-1}$. Then for every $C \in \U$ we will have that $\gon(C) = k$ and, in all
applications, $\U$ will have generic gonality $k$. We will find, in some cases with the help of of Proposition
\ref{mainbog}, a $(\U, L, k, \{A_C\}_{C \in \U})$-dominant set of divisors that, together with Proposition
\ref{lemma:somma}, will lead to the desired contradiction.

\begin{prop} 
\label{prop:finepiane}
Let $|L|$ be a base-point free complete linear system on an Enriques surface $S$ with $L^2 = 10$. Then a
general smooth curve $C \in |L|$ has $\gon(C) = 4$ and is not isomorphic to a smooth plane quintic.
\end{prop}

\begin{rem} {\rm Note that curves isomorphic to a smooth plane quintic do occur (on a proper closed
subset) by \cite{sta, ume} if $\phi(L) = 2$. On the other hand no smooth curve can be isomorphic
to a smooth plane curve of degree at least 6 by the results in \cite{kl2}}.
\end{rem}  

\begin{proof} Since $|L|$ is base-point free we have $2 \leq \phi(L) \leq 3$. Moreover, if $B > 0$ is a line
bundle such that $B^2 = 4$ then $B.L \geq 7$ by the Hodge index theorem, whence $\mu(L) \geq 5$.

As is well-known, for any smooth $C \in |L|$ we have $\gon(C) \leq \lfloor \frac{9}{2}
\rfloor = 4$. We will now prove that a general such curve cannot be trigonal. Suppose in fact that 
this is the case. By Lemmas \ref{cl:g=6,2} and \ref{cl:g=6,3} we see that we can apply Lemma
\ref{lemma:D} with $k = 3$, $M = E, N = L - E$, and we get the contradiction $3 \geq \min \{2\phi(L),
\mu(L) \} \geq 4$.

Now we show that $C$ is not isomorphic to a smooth plane quintic.

Suppose first that $\phi(L) = 2$. From Lemma \ref{cl:g=6,2} and 
Theorem \ref{cor:qnef}, we easily see that $|\O_C(E+E_1)|$ is a complete $g^1_5$ on $C$.  

If $C$ is isomorphic to a smooth plane quintic then it has no base-point free complete $g^1_5$. Hence
$|\O_C(E + E_1)|$ must have base points and it follows that $C$ passes through one of the two base points of
$|E + E_1|$. But this cannot happen for general $C \in |L|$.

Now suppose that $\phi(L) = 3$ and use Lemma \ref{cl:g=6,3}. 

By what we have just proved the open subset $\U' := |L|_{sm} - \Sigma_3$ is nonempty.
Consider now, in $|L|_{sm}$, the closed subset $\Sigma^2_5 = \{C \in |L|_{sm} : C \ \hbox{has a}\ 
g^2_5 \}$. We will prove that $\Sigma^2_5 \cap \U' \subset \U'$ is a strict inclusion. This will
give that on the nonempty open subset $\U' - \Sigma^2_5 \cap \U'$ every curve is not trigonal and is
not isomorphic to a plane quintic. 

Suppose that $\Sigma^2_5 \cap \U' = \U'$. Then every $C \in \U'$ has a $g^2_5$. Since $C$ has genus
$6$, we get that the $g^2_5$ must be very ample, therefore every $C \in \U'$ is isomorphic to a plane quintic.

Now certainly $\U'$ has generic gonality $4$. We now construct a particular $g^1_4$ on each $C \in \U'$. 

By Lemma \ref{cl:g=6,3}(i), $|\O_C(E_1 + E_2)|$ is a complete base-point free $g^2_7$. As $C$ cannot
be isomorphic to a smooth plane septic for reasons of genus, it follows that $\O_C(E_1 + E_2)$ is not very
ample. Since in addition any complete $g^1_5$ on $C$ must have a base point, there exists an effective divisor
$Z_3$ of degree $3$ on $C$ such that $|A_C:=\O_C(E_1 + E_2)(-Z_3)|$ is a complete $g^1_4$ on $C$, which has to
be base-point free, since $\gon(C) = 4$.

By Lemma \ref{cl:g=6,3}(i), it is easily seen that  $D$ is $(\U', L, 4, \{A_C\}_{C \in \U'})$-dominant. Now
$D^2 = 4$, $D.L = 7$ and $h^1(L-2D) \leq 1$ by Lemma \ref{cl:g=6,3}(ii).  By Proposition \ref{lemma:somma} we
see that equality occurs in (\ref{eq:databe}). But now Proposition \ref{lemma:somma}(c) yields a contradiction. 
\end{proof}

\begin{prop} 
\label{casi14,16}
Let $|L|$ be a base-point free complete linear system on an Enriques surface $S$ with $L^2 \geq 14$ and 
$\phi(L) = 2$. Then a general smooth curve $C \in |L|$ has gonality $4$ and possesses a
unique $g^1_4$, which is cut out by a unique genus one pencil on $S$. Moreover $W^1_4(C)$ is smooth 
when $L^2 = 14$. 
\end{prop}

\begin{rem} 
In fact $W^1_4(C)$ is smooth also when $L^2 = 16$. {\rm This will be proved in \cite{kl2}}.
\end{rem}

\begin{proof}
It is easy to see that there is a unique genus one pencil $|2E|$ such that
$E.L=\phi(L)$, thus also cutting out a $g^1_4$ on $C$. 

Now by \cite[Thm.1.4]{glmPN} we know that $\gon(C) > 3$, whence $\gon(C) = 4 = 2 \phi(L)$.
Moreover, again by \cite[Thm.1.4(1.7)]{glmPN}, we know that, if $L^2 \geq 18$, then every $g^1_4$ is
cut out by a genus one pencil on $S$ and we are done in this case.

Let us now treat the case $L^2 = 14$.  

Consider the nonempty open subset $\U = |L|_{sm} - \Sigma_3$ and suppose, to get a 
contradiction, that $\U$ has generic gonality $4$, that is that every $C \in \U$ has a $g^1_4 = 
|A_C|$ which is not cut out by a genus one pencil on $S$. 

Letting $L \sim 3E + E_1 + E_2$ be as in Lemma \ref{lemma:g=8,2}, we have that 
$|D| := |3E + E_1|$ is base-component free by \cite[Prop.3.1.6]{cd}. By Riemann-Roch $h^1(E_2 +
K_S) = 0$, whence, for every $C \in \U$, either $D_{|C} \geq A_C$ or $(E_2 + K_S)_{|C} \geq A_C$
by Lemma \ref{cor:nuovo}. In the second case we get a contradiction since $h^0((E_2 + K_S)_{|C}) =
h^0(E_2 + K_S) = 1$. As $|D|$ is base-component free, we have that $D$ is $(\U, L, k, \{A_C\}_{C \in
\U})$-dominant. But $\phi(D) = E.D = 1$, contradicting Proposition \ref{lemma:somma}(a). 

Therefore there is a curve $C_0 \in |L|_{sm}$ such that every $g^1_4$ on $C_0$ is cut out
by a genus one pencil on $S$, which is unique by what we said above. 

Let $A = (2E)_{|C_0}$ be the unique $g^1_4$ on $C_0$. We will prove that $\mu_{0, A}$ is
surjective. As is well-known \cite[Prop.IV.4.2]{acgh} this means that $W^1_4(C_0)$ is smooth at its
unique point $A$ and therefore the same holds in an open neighborhood of $C_0$ in $|L|$.

To see the surjectivity of $\mu_{0, A}$ we observe that, as $C_0$ is nontrigonal, if $\mu_{0, A}$ is
not surjective then, by the base-point-free pencil trick we have that $h^0(\omega_{C_0}-2A) \geq 3$. Since
$\deg(\omega_{C_0}-2A) = 6$ and $\Cliff(C_0) = 2$ we deduce that $|\omega_{C_0}-2A|$ is a base-point free
$g^2_6$ on $C_0$. Let $\varphi : C_0 \to X \subset \PP^2$ be the morphism defined by $|\omega_{C_0} - 2A|$. If
$\varphi$ is not birational, since $C_0$ is nontrigonal, then $\varphi$ is a cover of degree 2 of a smooth plane
cubic $X$, whence $C_0$ is bielliptic. But this is excluded since $C_0$ has a unique $g^1_4$.
Therefore $\varphi$ is birational, and again, since $C_0$ is nontrigonal, the image $X$ is a plane
sextic with two (possibly infinitely near) double points. Hence $C_0$, in its canonical embedding,
is isomorphic to a quadric section of a Del Pezzo surface, namely the anticanonical embedding
of the plane blown-up at two points. But in \cite[Lemma5.13]{klGM} it is proved that this case cannot
occur. 

Next we treat the case $L^2 = 16$. 

Consider the nonempty open subset $\U = |L|_{sm} - \Sigma_3$ and the finitely many
base-component free line bundles $L'$ such that $0 \leq L' \leq L$ and $\Bs |L'| \not= \emptyset$. We
shrink $\U$, if necessary, to a nonempty open subset $\U' \subseteq \U$ by removing the finitely
many closed subsets given by curves $C \in \U$ such that $C \cap \Bs |L'| \not= \emptyset$ for some
$L'$ as above.

Let $C \in \U'$, let $A_C$ be a $g^1_4$ on $C$ and apply Proposition \ref{mainbog} to
$A_C$. Then we must be in case (a) of the same proposition and $L \sim N_C +
N_C'$ with the three possibilities $(N_C')^2 = 0$ with $N_C'.L \leq 4$, or $(N_C')^2 = 2$ with
$N_C'.L \leq 6$, or $(N_C')^2 = 4$ with $N_C'.L = 8$.

In the first case $|N_C'|$ is a genus one pencil such that $A_C \sim \O_C(N_C')$. In the second case the Hodge
index theorem yields $N_C'.L = 6$. But $|N_C'|$ has two base points and by Proposition \ref{mainbog} we have
$\Bs |N_C'| \cap C \neq \emptyset$. In the third case the Hodge index theorem yields $L \eqv 2N_C'$ whence
$\phi(N_C') = 1$.  But now Proposition \ref{mainbog} implies that $C \cap \Bs |N_C'| \not= \emptyset$, a
contradiction.

Hence on any curve $C \in \U'$ there is a unique $g^1_4$, namely $A = (2E)_{|C}$.
\end{proof}

\section{The proof of Theorem \ref{main}} 
\label{sec:mainproof}

\begin{lemma}
\label{lemma:G}
Let $|L|$ be a base-component free complete linear system on an Enriques surface $S$ with $L^2 > 0$ and
let $C \in |L|_{sm}$. Then
$\gon(C) \leq \min\{2 \phi(L), \mu(L), \lfloor \frac{L^2}{4} \rfloor + 2\}$.
\end{lemma}

\begin{proof}
We can of course assume $L^2 \geq 4$ and that $k := \mu(L) \leq 2\phi(L) - 1$ and $L^2 \geq 4k - 4$. By
Proposition \ref{prop:mu} and Lemma \ref{lemmaprop} we have $k \geq 10$ and $L^2 \geq 36$. Pick $B > 0$ with
$B^2 = 4$ and $B.L = k + 2$. Then Riemann-Roch yields $h^0(L - B) \geq 9$. By Lemma \ref{lemma:A}(b) we have
$h^1(B) = h^1(L - B + K_S) = 0$, for otherwise we can find a line bundle $B_0 > 0$ such that $B_0^2 = 4$ and
$B_0.L \leq k$. Now if $\phi(B_0) = 1$ then, by Lemma \ref{lemma:nefred}, we can write $B_0 \sim 2F_1 + F_2$
with $F_i > 0, F_i^2 = 0$, $i = 1, 2$, whence $3 \phi(L) \leq B_0.L \leq k \leq 2\phi(L) - 1$, a contradiction.
Therefore $\phi(B_0) = 2$, but $B_0.L - 2 < k$, a contradiction.

Hence $h^0(B) = h^0(\O_C(B)) = 3$. Now $\O_C(B)$ cannot be very ample,
for reasons of genus and from the Hodge index theorem applied to $B$ and $L$. Therefore $\gon(C) \leq k =
\mu(L)$.
\end{proof}

\renewcommand{\proofname}{Proof of Theorem {\rm \ref{main}}}  
\begin{proof}

By Lemma \ref{lemma:G} we need to show that
\begin{equation} 
\label{eq:left}
\gon(C) \geq \min\{2 \phi(L), \mu(L), \lfloor \frac{L^2}{4} \rfloor + 2\}.
\end{equation}
holds for a general curve $C \in |L|$.

Now (\ref{eq:left}) holds if $\phi(L) = 1$. Therefore we can 
assume $\phi(L) \geq 2$ and $L^2 \geq 4$. If $4 \leq L^2 \leq 6$ the Hodge index
theorem shows that $\mu(L) \geq 3$, whence $\min\{2 \phi(L), \mu(L), \lfloor \frac{L^2}{4} \rfloor +
2\} = 3$ and a general curve in $|L|$ is trigonal by \cite[Prop.3.1.4, Prop.4.5.1,
Thm.4.5.4 and Rmk.4.5.2]{cd}. 

We can therefore assume that $L^2 \geq 8$ and that if $k$ is the gonality of a general curve 
in $|L|$ then $k \geq 3$. Moreover we can also assume that $L^2 \geq 4k - 4$ and $k \leq 2\phi(L) - 1$.
By Lemma \ref{lemma:D} we are left with the cases
\begin{equation}
\label{eq:I1}
L^2 = 4k - 2 \; \mbox{or} \; 4k - 4, \; \mbox{and} \ k \leq 2\phi(L) - 1.
\end{equation}
Hence $4k - 2 \geq L^2 \geq (\phi(L))^2 \geq \frac{(k + 1)^2}{4}$ gives $k \leq 13$ and $L^2 \leq 50$. By
Proposition \ref{poscone} we deduce that $k \leq 11$ and $L^2 \leq 42$.

We will prove \eqref{eq:left} by finding, for each case in \eqref{eq:I1}, an effective 
decomposition $L \sim M + N$ satisfying the conditions in Lemma \ref{lemma:D} or Lemma \ref{lemma:E} 
and then applying those lemmas. In each case we will first use \eqref{eq:I1} to determine $\phi(L)$ and then
Lemma \ref{lemma:nefred} to find a decomposition of $L$. Then we state what $M$ is, leaving most verifications
to the reader. 

\subsection{Notation} In the sequel of the proof we will let $E > 0$ be such that $E^2 = 0$ and 
$E.L = \phi(L)$. Moreover, any $E'$, $E_1$, $E_2$, etc., will be effective, nonzero isotropic
divisors.

\subsection{$k \leq 4$ and $L^2 = 4k - 2$} 
The possibilities are $(L^2, \phi(L), k) = (10, 2, 3)$, $(10, 3, 3)$ or $(14, 3, 4)$, 
and we are done by Lemma \ref{lemma:D} setting $M=E$.

\subsection{$k = 3$ and $L^2 = 8$}
We have $\phi(L) = 2$ and $L \sim 2E + E'$ with $E'.E = 2$, and one easily sees that either $\phi(L-E) = 2$
and $E'$ is primitive or $\phi(L-E) = 1$ and $E' \eqv 2E_1$. In the first case set $M=E$ and in the second $M
=E+E_1$.
 
\subsection{$k = 4$ and $L^2 = 12$} \label{ss:4,12} We have $\phi(L) = 3$. Choose a 
{\it maximal} $E$ such that $E.L = 3$. We have $L \sim 2E + E_1$ with $E.E_1 = 3$
and one easily sees that $E_1$ is primitive. We set $M=E$ and show how to verify that
$h^0(2N - L + K_S) =  h^0(E_1 + K_S) = 1$. By  Riemann-Roch and Theorem \ref{cor:qnef}
it suffices to show that $E_1$ is quasi-nef. 

Assume, to get a contradiction, that there is a $\Delta > 0$ such that $\Delta^2 = - 2$ and
$\Delta.E_1 \leq - 2$. Set $m = \Delta.E_1 \geq 2$. By Lemma \ref{A} there is a primitive $A > 0$ such
that $A^2 = 0$ and $E_1 \sim A + m \Delta$. We have $0 \leq L.\Delta = 2E.\Delta - m$ by the nefness of
$L$, whence $2E.\Delta \geq m$. From $3 = E.E_1 = E.A + m E.\Delta$ we get that $m = 2$,
$E.\Delta = E.A = 1$. Therefore we can write $L \sim 2(E + \Delta) + A$, contradicting the maximality
of $E$.

\subsection{$k = 5$ and $L^2 = 16$} 
We have $\phi(L) = 3$ or $4$. If $\phi(L) = 3$, let $L \sim 2E + E_1 + E_2$ be as 
in Lemma \ref{lemma:g=9,3}. If $h^0(2E + E_1 - E_2) = h^0(2E + E_1 - E_2 + K_S) = 1$ we set $M = E_2$. 
Otherwise, by Lemma \ref{lemma:g=9,3}, we must have that $E_1 \eqv 2F$ and we set 
$M = E + F + E_2$. (Note that $h^0(2M - L + K_S) = 1$ since 
$(2M - L + K_S).L = E_2.L < 2\phi(L)$.)

If $\phi(L) = 4$, by Lemma \ref{lemma:lat3} we have $L \eqv 2(E + E')$, with $E.E' = 2$
and we set $M = E + E'$. 

\subsection{$k = 5$ and $L^2 = 18$} 
We have $\phi(L) = 3$ or $4$. If $\phi(L) = 4$, set $M = E$.
If $\phi(L) = 3$ and $\phi(L - 2E) = 1$, then $L \sim 3(E + E_1)$ with $E.E_1=1$ and we set 
$M = 2E + E_1$. If $\phi(L) = 3$ and  $\phi(L - 2E) = 2$, then $L - 2E \sim E_1 + E_2 + E_3$, 
with $E_i.E_j = E.E_i=1$, for all $i \neq j$, and we set $M = E + E_1+ E_2$.   

\subsection{$k = 6$ and $L^2 = 20$}  
We have $\phi(L) = 4$ and it is easily seen that we can write $L \sim E_0 + E_1 + E_2 + E_3 + E_4$,
with $E_0 := E$ and $E_i.E_j = 1$ for $i \neq j$.
Therefore $\{E_0, \ldots , E_4\}$ is an isotropic 5-sequence and by \cite[Cor.2.5.6]{cd} we can find a 
divisor $F$ such that $F^2 = 0$, $F.E_i = 1$ for $0 \leq i \leq 4$. Then $F.L = 5$. Among all 
such $F$'s we choose a maximal one and let this be $M$. 
(We show how to see that $h^0(2N-L+K_S)=h^0(L+K_S-2M)=1$: Note that since $(L + K_S - 2M)^2 
= 0$ and $(L + K_S - 2M).L = 10$, we have that $L + K_S - 2M > 0$ and it is easily seen to be primitive.
Therefore, by Theorem \ref{cor:qnef}, it suffices to show that $L + K_S - 2M$ is quasi-nef. Assume 
there is a $\Delta > 0$ such that $\Delta^2 = -2$ and $\Delta.(L + K_S - 2M) \leq -2$. Then $\Delta.M
\geq 1$ by  the nefness of $L$ and if $\Delta.M = 1$ we get $\Delta.L = 0$ whence $(M + \Delta)^2 =
0$, $(M + \Delta).L = 5$, contradicting the maximality of $M$. Therefore $\Delta.M \geq 2$ and 
$(M + \Delta)^2 \geq 2$ whence $h^0(M + \Delta) \geq 2$ and we get $5 + \Delta.L = (M + \Delta).L \geq 
2\phi(L) = 8$, thus giving $\Delta.L \geq 3$. Now the nefness of $L$ implies that $\Delta.M \geq 3$, so 
that $(2M + \Delta)^2 \geq 10$. Also $(2M + \Delta).(L - 2M - \Delta) = 12 + \Delta.(L - 2M) - 
2 \Delta.M \leq 4$ and $(2M + \Delta).(L - 2M - \Delta) \geq 0$ by Lemma \ref{lemma10} since $L - 2M -
\Delta > 0$ by Lemma \ref{A}. Now the Hodge index theorem on $2M + \Delta$ and $L - 2M - \Delta$ yields a
contradiction.

\subsection{$k = 6$ and $L^2 = 22$} 
We have $\phi(L) = 4$ and it is easily seen that $L \sim 2E + E_1 + E_2 + E_3$ with $E.E_1 = 2$, $E.E_i = 1$
for $i = 2, 3$, $E_i.E_j = 1$ for $1 \leq i < j \leq 3$. Then we set $M = E_1$.

\subsection{$k = 7$ and $L^2 = 24$} 
We have $\phi(L) = 4$. One easily verifies that $(L-2E)^2=8$ and $\phi(L-2E)=2$, so that
$L - 2E \sim 2E_1 + E_2$ with $E_1$ primitive and $E_1.E_2 = 2$.

We have $4 = 2E.E_1 + E.E_2$. Hence, by Lemma \ref{lemma10}, we must have $E.E_1 = 1$ or $2$.
In the latter case we get $E.E_2 = 0$, then $E \eqv qE_2$ for some $q \geq 1$ by Lemma \ref{lemma10}.
From $E_1.E = E_1.E_2 = 2$ we get that $E \eqv E_2$ and we can set $M = 2E + E_1$.

If $E.E_1 = 1$ and $E.E_2 = 2$, setting $B = E + E_1$ and $A = E_2$, we have a
decomposition 
\begin{equation} 
\label{eq:7,24}
L \sim 2B + A, \: \mbox{with} \ A^2 = 0, \; B^2 = 2 \; \mbox{and} \  A.B = 4.
\end{equation}

Now among all effective decompositions of the form \eqref{eq:7,24}
pick one such that $B$ is maximal. As in case \ref{ss:4,12} one easily verifies that
$A$ is quasi-nef. If $h^1(A + K_S) > 0$, by Theorem \ref{cor:qnef}  and the fact that $A.L = 2\phi(L)$, we
have  that $A \eqv 2 A'$ for some divisor $A' > 0$ with $(A')^2 = 0$. In this case $L \eqv 2M$, setting 
$M = B + A'$. If $h^1(A + K_S) = 0$ we set $M = B + A$.

\subsection{$k = 7$ and $L^2 = 26$} 
By \eqref{eq:I1} and Proposition \ref{poscone} we get $\phi(L) = 4$. We have 
$L - 3E \sim E_1 + E_2$, with $E_1.E_2 = 1$. By symmetry we have the two
possibilities $(E.E_1, E.E_2) =  (2, 2)$ and $(1, 3)$. We set $M = 2E + E_2$.

\subsection{$k = 8$ and $L^2 = 28$} 
We have $\phi(L) = 5$. By Proposition \ref{poscone} and Lemma \ref{lemmaprop} we can
write $L \sim 2E + 2E_1 + E_2$ with $E.E_1 = E_1.E_2 = 2$ and $E.E_2 = 1$. We set $M = E + E_1 + E_2$.
(Note that $h^0(2M - L + K_S) = h^0(E_2 + K_S) = 1$, since $E_2.L = 6 < 2\phi(L)$.)
 
\subsection{$k = 8$ and $L^2 = 30$} 
We have $\phi(L) = 5$ and $L \sim 2E+ E_1 + E_2 + E_3$, with 
$E_1.E_2 = 1$, $E_1.E_3 = E_2.E_3 = 2$ and the three possibilities 
$(E.E_1, E.E_2, E.E_3) = (1, 1, 3)$, $(1, 2, 2)$ and $(2, 2, 1)$. In the first case we set $M = 2E + E_3$ and in
the two latter $M=E_2$.

\subsection{$k = 9$ and $L^2 = 32$} 
We have $\phi(L) = 5$ and $L \sim 3E+ E_1 + E_2$, with  $E.E_1=2$, $E.E_2= 3$ and $E_1.E_2 = 1$. We set $M = 2E
+ E_2$. (Note that $h^0(2M - L + K_S) \leq 1$ as $(2M - L + K_S).L = 8 < 2\phi(L)$.) 

\subsection{$k = 9$ and $L^2 = 34$} 
We have $\phi(L) = 5$ and one easily shows that $\phi(L - 2E) = 3$, so that
$L - 2E \sim 2E_1 + E_2 + E_3$, with $E_1.E_3 = 2$ and $E_1.E_2 = E_2.E_3 = 1$. 
Also either $E.E_1 = 2$ and $L \eqv 3E + 2E_1 + E_2$, in which case we set $M = 2E + E_1$ or
$E.E_1 = 1$, $(E.E_2, E.E_3) = (1, 2)$ or $(2, 1)$ and we set $M = E + E_1 + E_2$.

\subsection{$k = 10$ and $L^2 = 36$} 
We have $\phi(L) = 6$. By Lemmas \ref{lemma:lat3} and \ref{lemmaprop} we get $k = \mu(L) = 10$ and we are done.

\subsection{$k = 10$ and $L^2 = 38$} 
This  case cannot exist by \eqref{eq:I1} and Proposition \ref{poscone}. 

\subsection{$k = 11$ and $L^2 = 40$} 
We have $\phi(L) = 6$. By Proposition \ref{poscone} and Lemma \ref{lemmaprop} we have
that either $L \sim 3E + 2E_1 + E_2$ with $E.E_1 = E.E_2 = 2$, $E_1.E_2 = 1$ or 
$L \eqv 2D$ for a $D > 0$ with $D^2 = 10$. In the second case we set $M=D$. 
In the first case we set $M = 2E + 2E_1$. (Note that $h^0(2M - L + K_S) =  h^0(E + 2E_1 - E_2 + K_S) \leq 1$
since $(E + 2E_1 - E_2 + K_S).(E + E_1)  = 3 < 2\phi(E + E_1)$.) 

\subsection{$k = 11$ and $L^2 = 42$} 
We have $\phi(L) = 6$ and one can check that $\phi(L - 2E) = 4$. Then 
$L - 2E \sim 2E_1 + E_2 + E_3$, with $E_1.E_2 = E_1.E_3 = 2$ and $E_2.E_3 = 1$.

From $6 = E.L = 2E.E_1 + E.E_2 + E.E_3$ we find $E.E_2 \leq 3$ and $E.E_3 \leq 3$,
but if $E.E_i = 3$ for $i = 2$ or $3$, then $(E + E_i)^2 = 6$ yields the contradiction
$3\phi(L) = 18 \leq (E + E_i).L = 17$. Hence $E.E_2 \leq 2$ and $E.E_3 \leq 2$, so that $(E.E_1, E.E_2, E.E_3)
= (1, 2, 2)$ or $(2, 1, 1)$. We set $M = 2E + E_1 + E_2$.

This concludes the proof of \eqref{eq:left}, whence that of Theorem \ref{main}. 
\end{proof}
\renewcommand{\proofname}{Proof}

\renewcommand{\proofname}{Proof of Corollary {\rm \ref{main2}}}  
\begin{proof}
Apply Theorem \ref{main}, Propositions \ref{poscone} and \ref{prop:mu}, Lemmas
\ref{lemmaprop} and \ref{lemma:bassi2}.
\end{proof}
\renewcommand{\proofname}{Proof}

\section{The minimal gonality of a smooth curve in a complete linear system}
\label{sec:mingon}

In this section we prove Corollary \ref{vargon} and give some examples.

\renewcommand{\proofname}{Proof of Corollary {\rm \ref{vargon}}}                     

\begin{proof}
Let $k:= \gengon |L|$. Assume that there is a smooth curve $C_0 \in |L|$ with $k_0 := \gon(C_0) 
\leq k - 2$. Then $k_0 \leq \lfloor \frac{L^2}{4} \rfloor$ by Lemma \ref{lemma:G}, so that $L^2 \geq
4k_0 \geq 8$. By Remark \ref{usefulrem} with $b = k_0 - 1$ we get $L^2 \geq 2k_0 - 2 + 2 \phi(L)$. Fix
an ample divisor $H$ and let $A_0$ be a $g^1_{k_0}$ on $C_0$. By Proposition \ref{mainbog}(a) we obtain
an effective nontrivial decomposition $L \sim N + N'$ with $N.N' \leq k_0$. Moreover we have $(N -
N')^2 = L^2 - 4N.N' \geq 0$, and $H.(N - N') \geq 0$, so that by  Riemann-Roch we have that
either $N \eqv N'$ and $L \eqv 2N'$ or $L - 2N' = N - N' \geq 0$. Also recall  that $h^0(N') \geq 2$. Consider
the set
\[ X_L = \{D : \ h^0(D) \geq 2, D.(L - D) \leq k_0, (L - 2D)^2
\geq 0 \ \mbox{and either} \ L - 2D \geq 0 \ \mbox{or} \ L \eqv 2D \}. \]
Now $X_L \neq \emptyset$ since $N' \in X_L$, whence we can choose an element $D_0 \in X_L$ for which
$H.D_0$ is minimal. Note that $h^0(D_0) \geq 2$ implies that $D_0.L \geq 2\phi(L)$.

If $D_0^2 \geq 6$, then pick any $F > 0$ such that $F^2 = 0$ and $F.D_0 = \phi(D_0)$. Then $D_0 - F$
is easily seen to contradict  the minimality of $D_0$. Hence $D_0^2 \leq 4$. 

If $D_0^2 \leq 2$ we get
\[ k_0 \leq k - 2 \leq 2\phi(L) - 2 \leq D_0.L - 2 = D_0^2 + D_0.(L - D_0) - 2 \leq k_0 \]
whence $k_0 = k - 2$, $D_0^2 = 2$ and $D_0.L = 2\phi(L)$. Now the Hodge index theorem applied to
$D_0$ and $L$ implies that $\phi(L) \geq \lceil \sqrt{\frac{L^2}{2}} \rceil$.

If $D_0^2 = 4$ we first show that if $\phi(D_0) = 1$ then $L^2 = 8$ and $\phi(L) = 2$, 
whence $\phi(L) = \lceil \sqrt{\frac{L^2}{2}} \rceil$. If $\phi(D_0) = 1$ we can write $D_0
\sim 2F_1 +  F_2$ with $F_i > 0$, $F_i^2 = 0$, $i = 1, 2$, $F_1.F_2 = 1$ so that
\[ 3 \phi(L) \leq D_0.L = D_0^2 + D_0.(L - D_0) \leq k_0 + 4 \leq k + 2 \leq 2\phi(L) + 2 \]
whence $\phi(L) \leq 2$ and therefore $2 \leq k_0 \leq k - 2 \leq 2\phi(L) - 2 \leq 2$ giving $k = 4$, 
$k_0 = 2$ and $\phi(L) = 2$. Applying Proposition \ref{mainbog} to the hyperelliptic curve $C_0$ we get
that $L \sim N_0 + N'_0$  with $(N'_0)^2 \leq 2$ and $L.N'_0 \leq 4$. If $(N'_0)^2 = 0$ by
\cite[Lemma2.1]{glmPN} we have that $N'_0 \sim  2E$ is a genus one pencil such that $(2E)_{|C_0} \sim
A_0$, that is $E.L = 1$, a contradiction. Therefore $(N'_0)^2 = 2$ and the Hodge index theorem implies
that $L^2 \leq 8$.

Now if $D_0^2 = 4$ and $\phi(D_0) = 2$, we have, by Lemma \ref{lemma:G},
\[ k \leq \mu(L) \leq D_0.L - 2 = 2 + D_0.(L - D_0) \leq 2 + k_0 \leq k \]
whence $k_0 = k - 2$ and $\mu(L) = D_0.L - 2 = k \leq 2\phi(L)$. The Hodge index theorem applied to
$D_0$ and $L$ then gives $4 L^2 \leq (2\phi(L) + 2)^2$ and it is esily checked that 
$\phi(L) \geq \lceil \sqrt{\frac{L^2}{2}} \rceil$. 
\end{proof}

\subsection{A few examples}
\label{ex}
The flexibility of the Picard group of an Enriques surface allows us to give several examples to show 
the behavior of the gonality. We give two here, showing that all cases in Corollary \ref{vargon}, $\mingon |L|
= \gengon |L| - 2, \gengon |L| - 1, \gengon |L|$ do occur.  

\noindent {\bf Example 1}.
Let $|2E_i|$, $i = 1, 2$ be two genus one pencils on an Enriques surface $S$ such that 
$E_1.E_2 = 1$ and consider the line bundle 
\[ L \sim a E_1 + b E_2 \ \mbox{for} \ b \geq a \geq 3. \] 

It is easily checked that $2 \phi(L) = 2a \leq \min\{ \mu(L), \lfloor \frac{L^2}{4} \rfloor + 
2 \}$, so that, by Theorem \ref{main}, we have that $k:= \gengon |L| = 2a$. Moreover $2 \phi(L) <
\mu(L)$. 

We recall that by \cite[Prop.3.1.6, 3.1.4 and Thm.4.4.1]{cd} the base locus of $|E_1 + E_2|$
consists of two distinct points $x$ and $y$. Let $B \in |E_1 + E_2|$ be general and consider the exact
sequence
\[ 0 \hpil L - E_1 - E_2 \hpil \I_{\{x, y\}/S} \otimes L \hpil (L - E_1 - E_2)_{|B} \hpil 0, \]
the base point freeness of $L - E_1 - E_2 \sim (a - 1) E_1 + (b - 1) E_2$ (by 
\cite[Prop.3.1.6, 3.1.4 and Thm.4.4.1]{cd}) and the fact that $H^1(L - E_1 - E_2) = 0$, we get that
the base-scheme of $|\I_{\{x, y\}/S} \otimes L|$ is $\{x, y \}$. Therefore, by Bertini, 
the general curves $C_0 \in |\I_{\{x, y\}/S} \otimes L|$ and $C_0' \in
|\I_{\{x \}/S} \otimes L|$ are smooth and irreducible. 

Suppose that $b = a$. In this case we claim that both $k_0 = k - 2$ and $k_0 = k - 1$ do
occur as gonality of some smooth curve in $|L|$. Now $|(E_1 + E_2)_{|C_0}(- x -
y)|$ is a $g^1_{k - 2}$ on $C_0$, while $|(E_1 + E_2)_{|C_0'}(- x)|$ is a $g^1_{k - 1}$ on $C_0'$.
Therefore $k_0 := \gon(C_0) = k - 2$ by Corollary \ref{vargon}. 

To see that $k_0' := \gon(C_0') = k - 1$ we we first note that $y \not\in C_0'$ and also
that $\Bs|E_1 + E_2 + K_S| \cap C_0' = \emptyset$. By Corollary \ref{vargon} we can suppose that $k_0' =
k - 2$ and let $A$ be a $g^1_{k - 2}$ on $C_0'$. Applying Proposition \ref{mainbog}(a) we obtain an
effective nontrivial decomposition $L \sim N + N'$ with $|N'|$ base-component free, $N'.L \leq (N')^2 + 2a
- 2 \leq 4a - 4$ and $N'_{|C_0'} \geq A$.  

If $(N')^2 \geq 6$ we get $3a = 3 \phi(L) \leq N'.L \leq (N')^2 + 2a - 2$, so that $(N')^2
\geq a + 2$. But then the Hodge index theorem applied to $N'$ and $L$ gives a contradiction.

If $(N')^2 = 4$ we get that $L.N' \leq 2a + 2$, whence the contradiction $\mu(L) \leq 2a = 2
\phi(L)$. If $(N')^2 = 0, 2$ we have that $2a = 2 \phi(L) \leq L.N' \leq (N')^2 + 2a - 2$, whence
$(N')^2 = 2$ and $N'.L = 2a$, so that the Hodge index theorem implies that $L \eqv aN'$, whence $N' \eqv
E_1 + E_2$. On the other hand $h^0(N'_{|C_0'}) = h^0(N') = 2$, since $h^i(N' - L) = 0$ for $i = 0, 1$.
Therefore there are  two points $x', y' \in C_0'$ such that $N'_{|C_0'} \sim A + x' + y'$ and we deduce that
$x', y' \in \Bs|N'_{|C_0'}| = \Bs|N'| \cap C_0'$, contradicting our choice of $C_0'$. 

Now suppose that $b \geq a + 1$. In this case we claim that we cannot have $k_0 = k - 2$, 
while $k_0 = k - 1$ occurs as gonality of some smooth curve in $|L|$ if and only if $b = a + 1$. In fact,
as in the proof of Corollary \ref{vargon}, if $k_0 = k - 2$ occurs then there exists a divisor $D_0$ such that
$D_0^2 = 2$ and $D_0.L = k = 2a$. Now $D_0 \sim F_1 + F_2$ with $F_i > 0$, $F_i^2 = 0$, $i = 1, 2$, $F_1.F_2 =
1$, giving $2 \phi(L) \leq F_1.L + F_2.L = D_0.L = 2 \phi(L)$, whence the contradiction $F_1 \eqv E_2 \eqv
F_2$. Finally if $k_0 = k - 1$ occurs, then, as in the proof of Corollary \ref{vargon}, there exists a divisor
$D_0$ such that $D_0^2 = 2$ and $2a = k \leq D_0.L \leq k + 1 = 2a + 1$. Again $D_0 \sim F_1 + F_2$ with $F_i >
0$, $F_i^2 = 0$, $i = 1, 2$, $F_1.F_2 = 1$, giving $2 a \leq F_1.L + F_2.L = D_0.L \leq 2 a + 1$, which gives
necessarily $F_1 \eqv E_1$, $F_2 \eqv E_2$, whence $b = a + 1$. 

Then, when $b \geq a + 2$, minimal gonality and general gonality of $|L|$ coincide 
in this example, while, when $b = a + 1$, the case $k_0 = k - 1$ occurs since, as above, we can choose 
a smooth curve $C_0 \in |L|$ such that $x, y \in C_0$ and then $|(E_1 + E_2)_{|C_0}(- x - y)|$ is a
$g^1_{k - 1}$. 

\vskip .4cm

\noindent {\bf Example 2}.
Let $|2E_i|$, $i = 1, 2$, be two genus one pencils on an Enriques surface $S$ such that 
$E_1.E_2 = 2$ and consider the line bundle 
\[ L \sim a E_1 + a E_2 \ \mbox{for} \ a \geq 5. \] 

By Lemma \ref{lemmaprop} we have that $\mu(L) = 4a - 2 < \min\{ 2\phi(L), \lfloor \frac{L^2}{4} 
\rfloor + 2 \}$, whence $k:= \gengon |L| = 4a - 2$ by Theorem \ref{main}. Also $\phi(L) = 2a$.

We claim that both $k_0 = k - 2$ and $k_0 = k - 1$ occur as gonality of some smooth curve in
$|L|$. On any smooth curve $C \in |L|$ we have that $|(E_1 + E_2)_{|C}|$ is a $g^2_{k + 2}$ that cannot
be very ample, whence there are two points $x, y \in C$ such that $|(E_1 + E_2)_{|C}(- x - y)|$ is a $g^1_k$ on
$C$. If $C$ is general in $|L|$ this series computes the gonality of $C$ by Theorem \ref{main}. On the other
hand let $B_1, B_2 \in |E_1 + E_2|$ be two general smooth divisors and let $B_1 \cap B_2 = \{x_1, \ldots,
x_4\}$. As in Example 1, we find that the general curves $C_0 \in |\I_{\{x_1, \ldots,
x_4\}/S} \otimes L|$ and $C_0' \in |\I_{\{x_1, \ldots, x_3\}/S} \otimes L|$ are smooth and irreducible
and that $|(E_1 + E_2)_{|C_0}(- x_1 - \ldots - x_4)|$ is a $g^1_{k - 2}$ on
$C_0$, while $|(E_1 + E_2)_{|C_0'}(- x_1 - \ldots - x_3)|$ is a $g^1_{k - 1}$ on $C_0'$. Therefore the
minimal gonality of a smooth curve in $|L|$ is $k_0 := \gon(C_0) = k - 2$ by Corollary \ref{vargon}.
The fact that $k_0' := \gon(C_0') = k - 1$ can also be checked, with some cumbersome 
calculations.


\begin{thebibliography}{[ACGH]}

\bibitem[ACGH]{acgh} E.~Arbarello, M.~Cornalba, P.~A.~Griffiths, J.~Harris. 
\textit{Geometry of Algebraic Curves, Volume I}.
Grundlehren der Mathematischen Wissenschaften \textbf{267}.
Springer-Verlag, New York, 1985.

\bibitem[BEL]{bel} A.~Bertram, L.~Ein, R.~Lazarsfeld. 
\textit{Surjectivity of Gaussian maps for line bundles of large degree on curves}. 
Algebraic geometry (Chicago, IL, 1989), 15--25, LNM \textbf{1479}. Springer,
Berlin, 1991.

\bibitem[BPV]{bpv} W.~Barth, C.~Peters, A.~van de Ven. 
\textit{Compact complex surfaces}.
Ergebnisse der Mathematik und ihrer Grenzgebiete \textbf{4}. 
Springer-Verlag, Berlin-New York, 1984. 

\bibitem[Ca]{ca} G.~Castelnuovo. 
\textit{Sulle superfici algebriche le cui sezioni iperpiane sono curve iperellittiche}.
Rend. Circ. Mat. Palermo \textbf{4}, (1890) 73--88.

\bibitem[CD]{cd} F.~R.~Cossec, I.~V.~Dolgachev. 
\textit{Enriques Surfaces I}.
Progress in Mathematics \textbf{76}. Birkh\"auser Boston, MA, 1989.

\bibitem[Co1]{co1} F.~R.~Cossec. 
\textit{Projective models of Enriques surfaces}.
Math. Ann. \textbf{265}, (1983) 283--334.

\bibitem[Co2]{co2} F.~R.~Cossec. 
\textit{On the Picard group of Enriques surfaces}.
Math. Ann. \textbf{271}, (1985) 577--600.

\bibitem[CP]{cp} C.~Ciliberto, G.~Pareschi.
\textit{Pencils of minimal degree on curves on a $K3$ surface}. 
J. Reine Angew. Math. \textbf{460}, (1995) 15--36.

\bibitem[DM]{dm} R.~Donagi, D.~R.~Morrison. 
\textit{Linear systems on $K3$-sections}.
J. Diff. Geom. \textbf{29}, (1989) 49--64.

\bibitem[En]{en} F.~Enriques. 
\textit{Sui sistemi lineari di ipersuperficie algebriche ad intersezioni variabili iperellittiche}.
Math. Ann. \textbf{46}, (1895) 179--199.

\bibitem[F]{fr} R.~Friedman. 
\textit{Algebraic surfaces and holomorphic vector bundles}.
Universitext. Springer-Verlag, New York, 1998.

\bibitem[Fa]{fa} M.~L.~Fania. 
\textit{Trigonal hyperplane sections of projective surfaces}.
Manuscr. Math. \textbf{68}, (1990) 17--34. 

\bibitem[GL]{gl} M.~Green, R.~Lazarsfeld.
\textit{Special divisors on curves on a $K3$ surface}. 
Invent. Math. \textbf{89}, (1987) 357--370.

\bibitem[GLM]{glmPN} L.~Giraldo, A.~F.~Lopez, R.~Mu\~{n}oz.
\textit{On the projective normality of Enriques surfaces}.
Math. Ann. \textbf{324}, (2002) 135--158.

\bibitem[Ha]{ha} G.~M.~Hana.
\textit{Projective models of Enriques surfaces in scrolls.}
Math. Nachr. \textbf{279}, (2006) 242--254.

\bibitem[Ki]{kim} H.~Kim. 
\textit{Exceptional bundles on nodal Enriques surfaces}.
Manuscripta Math. \textbf{82}, (1994) 1--13.

\bibitem[KL1]{klvan} A.~L.~Knutsen, A.~F.~Lopez. 
\textit{A sharp vanishing theorem for line bundles on K3 or Enriques surfaces}. 
Proc. Amer. Math. Soc. \textbf{135}, (2007) 3495-3498.

\bibitem[KL2]{kl2} A.~L.~Knutsen, A.~F.~Lopez. 
\textit{Brill-Noether theory of curves on Enriques surfaces II}. 
In preparation.

\bibitem[KL3]{klGM} A.~L.~Knutsen, A.~F.~Lopez. 
\textit{Surjectivity of Gaussian maps for curves on Enriques surfaces}. 
Adv. Geom. \textbf{7}, (2007) 215--247.

\bibitem[KLM]{klm} A.~L.~Knutsen, A.~F.~Lopez, R.~Mu\~{n}oz. 
\textit{On the extendability of projective surfaces and a genus bound for Enriques-Fano
threefolds}.
Preprint math.AG/0605750.

\bibitem[Kn]{kn} A.~L.~Knutsen. 
\textit{On $k$th order embeddings of $K3$ surfaces and Enriques surfaces}.
Manuscr. Math. \textbf{104}, (2001) 211--237.

\bibitem[La]{la} R.~Lazarsfeld. 
\textit{Brill-Noether-Petri without degenerations}. 
J. Diff. Geom. \textbf{23}, (1986) 299--307.

\bibitem[Par]{par} G.~Pareschi. 
\textit{Exceptional linear systems on curves on Del Pezzo surfaces}.
Math. Ann. \textbf{291}, (1991) 17--38.

\bibitem[Pao]{pao} R.~Paoletti. 
\textit{Free pencils on divisors}.
Math. Ann. \textbf{303}, (1995) 109--123.

\bibitem[R]{re} M.~Reid. 
\textit{Special linear systems on curves lying on a K3 surface}. 
J. London Math. Soc. \textbf{13}, (1976) 454--458.

\bibitem[Se]{se} F.~Serrano. 
\textit{Extension of morphisms defined on a divisor}. 
Math. Ann. \textbf{277}, (1987) 395--413.

\bibitem[St]{sta} E.~Stagnaro. 
\textit{Constructing Enriques surfaces from quintics in $\PP^3_K$}. 
\textit{Algebraic Geometry - open problems (Ravello, 1982)}. LNM \textbf{997}.
Springer, Berlin, 1983.

\bibitem[SV]{sv} A.~J.~Sommese, A.~Van de Ven. 
\textit{On the adjunction mapping}. 
Math. Ann. \textbf{278}, (1987) 593--603.

\bibitem[T]{ty} A.~N.~Tyurin. 
\textit{Cycles, curves and vector bundles on an algebraic surface}. 
Duke Math. J. \textbf{54}, (1987) 1--26.

\bibitem[Um]{ume} Y.~Umezu. 
\textit{Normal quintic surfaces which are birationally Enriques surfaces}. 
Publ. Res. Inst. Math. Sci. \textbf{33}, (1997) 359--384.

\bibitem[Wa]{wa} J.~Wahl. 
\textit{Introduction to Gaussian maps on an algebraic curve}.
Complex Projective Geometry, Trieste-Bergen 1989, London Math. Soc. Lecture Notes Ser. 
\textbf{179}. CUP, Cambridge 1992, 304-323.  

\bibitem[Za]{za} F.~L.~Zak. 
\textit{Some properties of dual varieties and their applications in projective geometry}.
Algebraic geometry (Chicago, IL, 1989), 273--280. LNM \textbf{479}. Springer,
Berlin, 1991.

\end{thebibliography}
\end{document}